\documentclass{gtpart}   
\theoremstyle{plain}
\usepackage{amsmath, amsfonts, amsthm, amssymb,amscd}
\usepackage{graphics}
\usepackage{color}
\usepackage{epsfig}
\usepackage[all]{xy}
\graphicspath{ {Figures/} }

\title[Khovanov and Heegaard Floer homology]{Khovanov homology, sutured {F}loer homology, and annular links}
\author{J. Elisenda Grigsby}
\address{Boston College Mathematics Department\\301 Carney Hall\\Chestnut Hill, MA 02467}
\email{grigsbyj@bc.edu}

\author{Stephan M. Wehrli}
\address{Syracuse University Mathematics Department\\215 Carnegie\\Syracuse, NY 13244}
\email{smwehrli@syr.edu}


\arxivreference{math.GT/0907.4375}  %%% please supply if 
\arxivpassword{u9ut4}   %%% paper is in the arXiv

\theoremstyle{plain}
\newtheorem{theorem}{Theorem}[section]
\newtheorem{lemma}[theorem]{Lemma}
\newtheorem{proposition}[theorem]{Proposition}

\newtheorem*{theoremReinterp}{Theorem~\ref{thm:Reinterp}}
\newtheorem*{theoremCutting}{Theorem~\ref{thm:Cutting}}
\theoremstyle{definition}
\newtheorem{notation}[theorem]{Notation}
\newtheorem{definition}[theorem]{Definition}
\newtheorem{remark}[theorem]{Remark}
\newtheorem{example}[theorem]{Example}
\newcommand{\Ozsvath}{{Ozsv{\'a}th} }

\newcommand{\Juhasz}{{Juh{\'a}sz} }

\newcommand{\Qrat}{\ensuremath{\mathbb{Q}}}
\newcommand{\Real}{\ensuremath{\mathbb{R}}}
\newcommand{\Zint}{\ensuremath{\mathbb{Z}}}

\newcommand{\bL}{\ensuremath{\mathbb{L}}}

\newcommand{\OO}{\ensuremath{\mathbb{O}}}
\newcommand{\XX}{\ensuremath{\mathbb{X}}}

\newcommand{\boldalpha}{\ensuremath{\mbox{\boldmath $\alpha$}}}
\newcommand{\boldbeta}{\ensuremath{\mbox{\boldmath $\beta$}}}

\newcommand{\boldSigma}{\ensuremath{\mbox{\boldmath $\Sigma$}}}
\newcommand{\cP}{\ensuremath{\mathcal{P}}}
\newcommand{\cI}{\ensuremath{\mathcal{I}}}
\newcommand{\cV}{\ensuremath{\mathcal{V}}}
\newcommand{\cF}{\ensuremath{\mathcal{F}}}

\newcommand{\cM}{\ensuremath{\mathcal{M}}}
\newcommand{\Ztwo}{\ensuremath{\mathbb{Z}_2}}
\newcommand{\Wedge}{\ensuremath{\Lambda}}
\newcommand{\Torus}{\ensuremath{\mathbb{T}}}

\newcommand{\hatSigma}{\ensuremath{\widehat{\Sigma}}}

\newcommand{\ZH}{\ensuremath{Z_{H}}}
\newcommand{\VbarH}{\ensuremath{\overline{V}_{H}}}
\newcommand{\VH}{\ensuremath{V_{H}}}
\newcommand{\hatbL}{\ensuremath{\widehat{\bL}}}
\newcommand{\Thetahat}{\ensuremath{\widehat{\Theta}}}

\begin{document}
\bibliographystyle{plain}

\begin{abstract} %In \cite{MR2141852} \Ozsvath and Szab{\'o} show how to associate to a link $\bL$ in $S^3$ a spectral sequence whose $E^2$ term is the reduced Khovanov homology of (the mirror of) $\bL$ and whose $E^\infty$ term is the Heegaard Floer homology of the double-branched cover of $\bL$.  
In \cite{GT07060741}, Lawrence Roberts, extending the work of \Ozsvath and Szab{\'o} in \cite{MR2141852}, showed how to associate to a link, $\bL$, in the complement of a fixed unknot, $B \subset S^3$, a spectral sequence whose $E^2$ term is the Khovanov homology of a link in a thickened annulus defined in \cite{MR2113902}, and whose $E^\infty$ term is the knot Floer homology of the preimage of $B$ inside the double-branched cover of $\bL$.  

In \cite{GT08071432}, we extended \cite{MR2141852} in a different direction, constructing for each knot $K \subset S^3$ and each $n \in \Zint_+$, a spectral sequence from Khovanov's categorification of the reduced, $n$--colored Jones polynomial to the sutured Floer homology of a reduced $n$--cable of $K$.  In the present work, we reinterpret Roberts' result in the language of Juh{\'a}sz's sutured Floer homology \cite{MR2253454} and show that the spectral sequence of \cite{GT08071432} is a direct summand of the spectral sequence of \cite{GT07060741}.
\end{abstract}
\maketitle

\section{Introduction}
Heegaard Floer homology \cite{MR2113019} and Khovanov homology \cite{MR1740682} have transformed the landscape of low-dimensional topology in the past decade, generating a wealth of applications, most notably to questions in knot concordance (cf., \cite{MR2026543}, \cite{GT0306378}, \cite{GT0402131}), Dehn surgery (cf., \cite{MR2168576}, \cite{GT08071341}), and contact geometry (cf., \cite{MR2153455}, \cite{GT0412184}).  The philosophies underlying the theories' constructions are quite different, yet there are intriguing connections between the two.  The best-understood such connection is an algebraic relationship, discovered by Ozsv{\'a}th-Szab{\'o}, between the homology of a Khovanov-type chain complex associated to a link and the homology of a Heegaard Floer-type chain complex associated to its double-branched cover.  

Specifically, in \cite{MR2141852} \Ozsvath and Szab{\'o} show how to associate to a link $\bL \subset S^3$ a spectral sequence whose $E^2$ term is $\widetilde{Kh}(\overline{\bL})$ and whose $E^\infty$ term is $\widehat{HF}(\boldSigma(S^3, \bL))$.   Here (and throughout), $\widetilde{Kh}$ denotes the reduced version of Khovanov homology defined in \cite{MR2034399}, $\overline{\bL}$ denotes the mirror of $\bL$, $\boldSigma(A, B)$ denotes the double-branched cover of $A$ branched over $B$, and $\widehat{HF}$ denotes (the hat version of) Heegaard Floer homology.  Unless explicitly stated otherwise, all Khovanov and Heegaard Floer homology theories discussed in this paper will be considered with coefficients in $\Zint_2$.  

%is induced by a filtered chain complex obtained by counting holomorphic polygons in a Heegaard multi-diagram associated to (a diagram for) $L$.  

%At the time of its original discovery, it was unclear both how general and how useful the spectral sequence connecting Khovanov and Heegaard Floer homology would be.  
Later work of Roberts, building on work of Plamenevskaya, demonstrated that this relationship was both useful and more general than originally believed.  Specifically, given a link $\bL$ in the complement of a fixed unknot, $B \subset S^3$, Roberts exhibits, in \cite{GT07060741}, the existence of a spectral sequence from $Kh^*(\overline{\bL})$ (where $Kh^*$ is a version of Khovanov homology for links in product manifolds defined in \cite{MR2113902}) to (a variant of) the knot Floer homology of $\widetilde{B} \subset \boldSigma(S^3, L)$, where $\widetilde{B}$ is the preimage of $B$ in $\boldSigma(S^3, \bL)$.  This allowed him to establish a relationship, first conjectured in \cite{GT0412184}, between Plamenevskaya's transverse invariant \cite{GT0412184} and Ozsv{\'a}th-Szab{\'o}'s contact invariant \cite{MR2153455}.  Baldwin and Plamenevskaya, in \cite{GT08082336}, used (an extension of) this relationship to establish the tightness of a number of non Stein-fillable contact structures.
%knot Floer homology of certain fibered knots satisfy nice properties.  In particular, he 
%This spectral sequence is a specific instance of a {\em link surgeries spectral sequence} (an iterated application of the long exact sequence associated to a knot), applied to the link obtained as the double-branched cover of a neighborhood of all of the knots obtai

In \cite{GT08071432}, we extended \cite{MR2141852} in a different direction, using a version of Heegaard Floer homology for sutured manifolds (Definition \ref{defn:suturedman}) developed by \Juhasz in \cite{MR2253454}, yielding a proof that a simple variant of Khovanov homology, categorifying the reduced, $n$--colored Jones polynomial and defined in \cite{MR2124557}, detects the unknot whenever $n \geq 2$.  

The goal of the present work is to propose a single general framework unifying all of these results.  As in \cite{GT08071432}, this general framework uses Gabai's sutured manifold theory \cite{MR723813} and Juh{\'a}sz's sutured Floer homology \cite{MR2253454}.  It not only incorporates all known results of this type but also clarifies their relationship to each other.  In particular, it can be shown to satisfy nice naturality properties with respect to certain TQFT operations, as detailed in \cite{SurfDecomp}.%, an invariant of balanced sutured manifolds (Definitions \ref{defn:suturedman}, \ref{defn:balancedSM}).  

More specifically, let $F$ be an oriented surface with $\partial F \neq \emptyset$, $F \times I$ a product sutured manifold (see Definition \ref{defn:suturedman} and Example \ref{example:product}), and $T \subset F \times I$ a tangle (properly imbedded $1$--manifold) that is both {\em admissible}, i.e., \[T \cap (\partial F \times I) = \emptyset,\] and {\em balanced}, i.e., \[|T \cap (F \times \{1\})| = |T \cap (F \times \{0\})|= n \in \Zint_{\geq 0}.\]

Then there should exist a spectral sequence whose $E^2$ term is an appropriate version of the Khovanov homology of $\overline{T}$, which we will denote $Kh^*(\overline{T})$, and whose $E^\infty$ term is the sutured Floer homology of $\boldSigma(F \times I, T)$, which we will denote $SFH(\boldSigma(F \times I, T))$.  We shall explore these ideas more fully in the sequel--in particular, the appropriate version of Khovanov homology for balanced tangles in product sutured manifolds should be similar to what appears in \cite{QA0410238}, with abelianized gradings--focusing here on the case when $F=A$ is an annulus and $T = \bL$ is a $0$--balanced tangle (link):

%We conjecture:

%\begin{conjecture}\label{conj:GenFramework}
%Let $F$ be an oriented surface with $\partial F \neq \emptyset$, $F \times I$ the associated product sutured manifold, and $T \subset F \times I$ an admissible, balanced tangle.  Then there exists a spectral sequence whose $E^2$ term is $Kh^*(\overline{T})$, an appropriate version of the Khovanov homology of $\overline{T}$, and whose $E^\infty$ term is the sutured Floer homology of $\boldSigma(F \times I, T)$.

%Furthermore, if $\gamma \subset F$ is a properly-imbedded, oriented arc representing a non-zero element of $H_1(F, \partial F)$, and $F'$ is the surface with boundary obtained by cutting $F$ along $\gamma$, and $T' \subset F' \times I$ is the balanced tangle obtained by performing a surface decomposition along $\gamma \times I$, then the spectral sequence associated to $T'$ is a direct summand of the spectral sequence associated to $T$.
%\end{conjecture}

\begin{theoremReinterp}
Let $\bL \subset A \times I$ be a link in the product sutured manifold $A \times I$.  Then there is a spectral sequence whose $E^2$ term is $Kh^*(\overline{\bL})$ and whose $E^\infty$ term is $SFH(\boldSigma(A \times I, \bL))$.
\end{theoremReinterp}

Note that we treated the case when $F=D$ in  \cite{GT08071432}:%\footnote{In \cite{GT08071432}, a Khovanov-type chain complex is constructed using projections of $T \subset D \times I$ to a vertical disk.  A straight-line isotopy moving the endpoints of $T$ to $\partial D_\pm$ yields an equivalent projection of $T$ to a horizontal disk with endpoints, $T \cap D_+$, marked ``+'' and endpoints, $T \cap D_-$, marked ``-''.  A horizontal projection of this type has the property that all ``+'' endpoints (similarly, all ``-'' endpoints) are adjacent  on $\partial D$.  More general horizontal projections can be considered.}

\begin{theorem} \label{thm:DiskThm} \cite[Prop. 5.20]{GT08071432} Let $T \subset D \times I$ be an admissible, balanced tangle.  Then there is a spectral sequence whose $E^2$ term is $Kh^*(\overline{T})$ and whose $E^\infty$ term is $SFH(\boldSigma(D \times I, T))$.\footnote{$Kh^*(\overline{T})$ is denoted $V(T)$ in \cite[Sec. 5]{GT08071432}.}
\end{theorem}

%Furthermore, observing that

%\begin{eqnarray*}
%  Kh^*(\overline{T}) &\cong& \widetilde{Kh}(L), \,\,\,\mbox{ and}\\ 
%  SFH(\boldSigma(D \times I, T)) &\cong& \widehat{HF}(\boldSigma(S^3, L))\\
%\end{eqnarray*}

%when $T$ is the $1$--balanced tangle obtained by removing the neighborhood of a point on $L$, we recover Ozsv{\'a}th-Szab{\'o}'s original result.

In fact, Theorem \ref{thm:Reinterp} is a reinterpretation (and modest extension) of Roberts' main result in \cite{GT07060741}.  (Note that in \cite{GT07060741}, Roberts restricts to those $\bL$ intersecting a spanning disk of $B$ in an odd number of points, while Theorem \ref{thm:Reinterp} requires no such restriction.)  To understand the connection, let $A$ be an oriented annulus, $I = [0,1]$ the oriented closed unit interval, and $\bL \subset A \times I$ a link, where $A \times I$ has been identified as the standard sutured complement of a standardly-imbedded unknot, $B \subset S^3$, via the identification: 
\begin{eqnarray*}
A \times I &=& \{(r,\theta, z) \,|\, r \in [1,2], \theta \in [0,2\pi), z \in [0,1]\} \subset \Real^3 \cup \infty = S^3,\\
B &=& \{(r,\theta, z)\,|\,r=0\} \cup \infty \subset S^3.%Note that as an oriented knot, identified with $\lambda \times I$, $B$ is oriented in terms of decreasing 
\end{eqnarray*}

In \cite[Prop. 1.1]{GT07060741}, Roberts constructs a spectral sequence from $Kh^*(\overline{\bL})$, where $Kh^*$ is a version of Khovanov homology for links in product manifolds defined in \cite{MR2113902}\footnote{$Kh^*(\bL)$ is denoted $H(\bL)$ in \cite[Sec. 2]{GT07060741}.} to (a variant of) the knot Floer homology of $\widetilde{B} \subset \boldSigma(S^3, \bL)$, where $\widetilde{B}$ is the preimage of $B$ in $\boldSigma(S^3, \bL)$.  We show (see Proposition \ref{prop:SFHFiltCpx}) that this (variant of) the knot Floer homology of $\widetilde{B}$ is, in fact, just the sutured Floer homology of $\boldSigma(A \times I, \bL)$.

%Roberts' result (compare \cite[Prop. 1.1]{GT07060741}) is restated as follows, in the language of sutured manifolds:

%\begin{theoremReinterp}
%Let $\bL \subset A \times I$ be a link in the product sutured manifold $A \times I$.  Then there exists a spectral sequence from $Kh^*(\overline{\bL})$ to $SFH(\boldSigma(A \times I, \bL))$.
%\end{theoremReinterp}

\begin{figure}
\begin{center}
\resizebox{4.5in}{!}{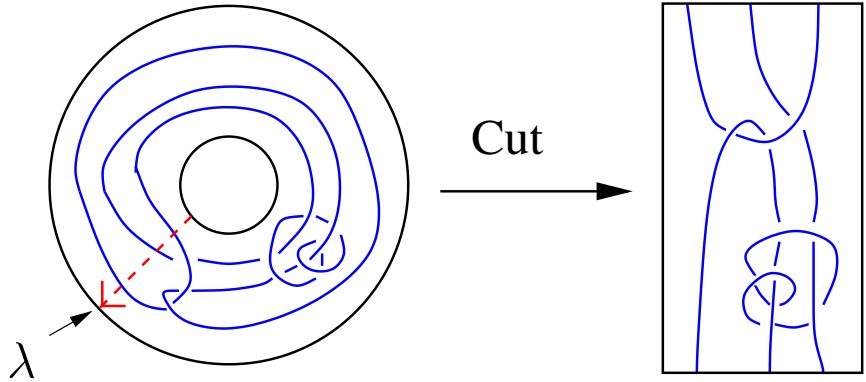}
\end{center}
\caption{Cutting an annular link projection along a homologically non-trivial oriented arc }
\label{fig:CutGlue}
\end{figure}

Furthermore, there is a nice relationship between the spectral sequences of Theorems \ref{thm:DiskThm} and \ref{thm:Reinterp}.  Specifically a link $\bL \subset A \times I$ can be cut along a vertical disk as in Figure \ref{fig:CutGlue} to form an admissible, balanced tangle $T \subset D \times I$.  The following is a corollary of \cite[Thm. 4.5]{SurfDecomp}.  

%if $D$ represents an oriented disk, then %there is a natural operation on (projections of) links in $A \times I$ producing (projections of) tangles in $D \times I$:  Specifically, 
%any projection of a link in $A \times I$ can be cut along a vertical disk of the form $\lambda \times I$, where $\lambda$ is a homologically non-trivial, properly-imbedded arc in $A$ transverse to the link projection as in Figure \ref{fig:CutGlue} to produce a tangle in $D \times I$.  Any tangle produced by this procedure is admissible and balanced.  An immediate corollary of \cite[Thm. 4.5]{GT08071432} is then:%}, i.e., \[T \cap (\partial D \times I) = \emptyset,\] and {\em balanced}, i.e., \[\mbox{i.e., } |T \cap (D \times \{1\})| = |T \cap (D \times \{0\})|= n \in \Zint_{\geq 0}.\]  In \cite{GT08071432}, we show how to associate to each (ambient isotopy class of) admissible balanced tangle $T \subset D \times I$ (Definition \ref{defn:balancedtangle}) a spectral sequence whose $E^2$ term is $Kh^*(\overline{T})$, an appropriate version of the Khovanov homology of $T$, and whose $E^\infty$ term is $SFH(\boldSigma(D \times I, T))$.  This spectral sequence is related to the one described in Theorem \ref{thm:Reinterp} as follows:

%Combined with the main result of \cite{SurfDecomp}, this will immediately imply that the spectral sequence we construct in \cite{GT08071432} can be understood as a direct summand of the spectral sequence Roberts constructs in \cite{GT08082817}.

\begin{theoremCutting} Let $\bL \subset A \times I$ be an isotopy class representative of an annular link admitting a projection, $\cP(\bL)$, and let $\lambda \subset A$ be a properly imbedded oriented arc representing a nontrivial element of $H_1(A,\partial A)$ such that $\lambda$ intersects $\cP(\bL)$ transversely.   % $\lambda \times I$ intersects $\bL$ transversely.  
Let $T \subset D \times I$ be the balanced tangle in $D \times I$ obtained by decomposing $A \times I$ (Definition \ref{defn:DecResult}) along the surface $\lambda \times I$ endowed with the product orientation.

Then the spectral sequence \[Kh^*(\overline{T}) \rightarrow SFH(\boldSigma(D \times I, T))\] is a direct summand of the spectral sequence \[Kh^*(\overline{\bL}) \rightarrow SFH(\boldSigma(A \times I, \bL)).\]

Furthermore, the direct summand is trivial if there exists some $\bL' \subset A \times I$ isotopic to $\bL$ satisfying \[|(\lambda \times I) \pitchfork \bL'| \lneq |(\lambda \times I) \pitchfork \bL|\]

%Furthermore, if there exists some $\bL' \subset A \times I$ isotopic to $\bL$ satisfying \[|(\lambda \times I) \pitchfork \bL'| \lneq |(\lambda \times I) \pitchfork \bL|,\] then the direct summand is trivial.
\end{theoremCutting}

%where two admissible $n$--balanced tangles are considered equivalent if they are ambiently isotopic through admissible $n$--balanced tangles.  

%if we view $D \times I$ and $A \times I$ as product sutured manifolds (Definition \ref{defn:suturedman}), this cutting procedure can be described in terms of a surface decomposition.

Theorem \ref{thm:Cutting} can be viewed as the first example of the ``naturality'' of the relationship between Khovanov homology and Heegaard Floer homology; under various natural geometric operations, the spectral sequence behaves ``as expected.''  See \cite{SurfDecomp} for more examples.

It is also interesting to note that, given a link, $\bL \subset S^3$, any unknot, $B \subset S^3 - N(\bL)$, endows the Khovanov chain complex associated to $\bL \subset B^3$ with a \Zint--filtration, via the identification \[S^3 - N(B) \leftrightarrow A \times I\] (see \cite{GT07060741} and Proposition \ref{prop:SFHFiltCpx}).
The extra grading inducing this $\Zint$--filtration has the following representation-theoretic interpretation. Suppose $T \subset D\times I$ is an $n$--balanced tangle obtained by decomposing $\bL \subset A \times I$ along $\lambda \times I$, as in Theorem~\ref{thm:Cutting}. Then the $U_q(\mathfrak{s}l_2)$ tangle invariant associated to $T$ is an endomorphism, $J(T)$, of the $U_q(\mathfrak{s}l_2)$-module $V^{\otimes n}$, where $V$ denotes the standard $2$--dimensional irreducible representation of $U_q(\mathfrak{s}l_2)$. Denoting by $a_k(T)\in\Zint[q^{\pm 1}]$ the trace of the restriction of $J(T)$ to the $k$-th weight space in $V^{\otimes n}$, one can show that the sum, \[\sum_{k\in Z}(tq)^ka_k(T)\in\Zint[q^{\pm 1},t^{\pm 1}],\] is precisely the graded Euler characteristic, $\chi_{t,q}(Kh^*(\bL))),$ of $Kh^*(\bL)$, defined as in \cite{GT07060741}, with the power of $t$ corresponding to the extra grading. Thus, the extra grading is related to the decomposition of $J(T)$ according to weight spaces in $V^{\otimes n}$.

%and the standard $2$--dimensional irreducible representation $V$ of $U_q(\mathfrak{s}l_2)$ is a $\Complex(q)$--linear map $J(T)\colon V^{\otimes n}\rightarrow V^{\otimes n}$ which intertwines the $U_q(\mathfrak{s}l_2)$ action and thus preserves the decomposition of $V^{\otimes n}$ into weight spaces.

%The extra grading inducing this $\Zint$--filtration should have an interesting representation-theoretic interpretation.  For example, in the case where $K \subset S^3$ is a knot, $\bL \subset S^3$ its $0$--framed $n$--cable, and $B$ an oriented meridian enclosing all strands of the cable, the grading induced by $B$ is related to the dimensions of the various irreducible subrepresentations of $U_q(\mathfrak{s}l_2)$ in $V^{\otimes n}$, where $V$ is the standard $2$--dimensional irreducible representation of $U_q(\mathfrak{s}l_2)$.

We also remark that choosing an alternative rational slope on $\partial (S^3 - N(B))$ corresponds to choosing an alternative Dehn filling of $B$, furnishing a plausible means of constructing Khovanov-type invariants for links in lens spaces.

{\flushleft {\bf Acknowledgments:}} This work has benefited greatly from conversations with Ken Baker, John Baldwin, Matt Hedden, Andr{\'a}s Juh{\'a}sz, Mikhail Khovanov, Rob Kirby, Robert Lipshitz, Peter Ozsv{\'a}th, Lawrence Roberts, and Liam Watson.  In addition, a large portion of this work was completed while both authors were visiting postdoctoral fellows at Columbia University.  We are grateful to the Columbia mathematics department for its hospitality.  The first author was partially supported by an NSF postdoctoral fellowship and NSF grant number DMS-0905848.  The second author was partially supported by a Swiss NSF fellowship for prospective researchers and a postdoctoral fellowship of the Fondation Sciences Math\'ematiques de Paris.

\section{Spectral Sequence from Khovanov to Sutured Floer}

%The purpose of this section is to prove Theorem \ref{thm:}, which asserts the existence of a spectral sequence whose $E^2$ term is the Khovanov-type homology associated to a link $\bL \subset A \times I$ and whose $E^\infty$ term is the sutured Floer homology of the double-branched cover, $\boldSigma(A \times I, \bL)$, of $A \times I$ over $\bL$.

%Since \cite[Prop. 4.1]{GT08071432} proves the existence of a spectral sequence coming from a certain natural link $L \subset \boldSigma(A \times I, \bL)$, our strategy (following \cite{MR2141852} and \cite{GT08071432}) will be to relate the $E^2$ term of this spectral sequence to the Khovanov homology of $\bL \subset A \times I$ by proving the equivalence of a certain Khovanov-type and sutured Floer-type ``functor'' on the category of resolved links in $A \times I$.

%In fact, Section \ref{sec:} will describe a splitting of this spectral sequence in terms of the evaluation of first Chern classes of relative Spin$^c$ structures on a natural imbedded surface.  Section \ref{sec:} will then describe the effect on the spectral sequence of decomposing along this surface.

This section is devoted to proving the main theorem:

\begin{theorem} \label{thm:Reinterp}
Let $\bL \subset A \times I$ be a link in the product sutured manifold $A \times I$.  Then there exists a spectral sequence from $Kh^*(\overline{\bL})$ to $SFH(\boldSigma(A \times I, \bL))$.
\end{theorem}

In the above, $Kh^*(\overline{\bL})$ is the homology of a Khovanov-type chain complex associated to $\overline{\bL}$, the mirror of $\bL$, %constructed in Section \ref{sec:KBSM}, 
while $SFH(\boldSigma(A \times I, \bL))$ is the {\em sutured Floer homology} of $\boldSigma(A \times I, \bL)$, the homology of a Heegaard Floer-type chain complex.  %constructed in Section \ref{sec:SFH}.  
Both constructions begin by associating a cube of resolutions to a projection (diagram) of $\bL \subset A \times I$.  As in \cite{GT08071432}, the main step in the proof will be a verification of the equivalence of Khovanov and sutured Floer ``functors'' applied to the vertices and edges of this cube.  We begin by fixing notation and recalling some standard definitions in Section \ref{sec:not}.  We go on to describe the Khovanov and sutured Floer functors on resolved link projections in Sections \ref{sec:KBSM} and \ref{sec:SFH}, establishing their equivalence in Proposition \ref{prop:Equiv}.

\subsection{Notation and Standard Definitions}\label{sec:not}
\begin{notation}Throughout, $A$ denotes an oriented annulus and $I := [0,1]$ denotes the oriented closed unit interval.

Whenever we write $A \times I$, we shall always assume we have fixed an identification: \[A \times I := \{(r,\theta,z)\,\,|\,\,r \in [1,2], z \in [0,1]\} \subset \Real^2 \times \Real.\}\]  Let $A_+$ (resp., $A_-$) denote $A \times \{1\}$ (resp., $A \times \{0\}$).  %More generally, let $A_t$ denote $A \times \{t\}$ for each $t \in [0,1]$. 
\end{notation}

Most of the following definitions can be found in \cite{MR723813}, \cite{MR2253454}, \cite{MR2390347}, and \cite{GT08071432}.  See also \cite{MR1740682} and \cite{MR2113019}.
%This theorem is a reinterpretation (and modest generalization) of \cite[]{}.  Our main tool will be a version of Ozsv{\'a}th-Szab{\'o}'s link surgeries spectral sequence for sutured manifolds, proved in \cite{GT08071432}.

%, along with a version of Juh{\'a}sz's surface decomposition theorem in the complement of a framed link, proved in \cite{SurfDecomp}.

%\ref{thm:Cutting}
%\subsection{Definitions and Notations} 

\begin{definition} \cite{MR723813} \label{defn:suturedman}
A {\em sutured manifold} $(Y,\Gamma)$ is a compact, oriented $3$--manifold with boundary $\partial Y$ along with a set $\Gamma \subset \partial Y$ of pairwise disjoint annuli $A(\Gamma)$ and tori $T(\Gamma)$.  The interior of each component of $A(\Gamma)$ contains a {\em suture}, an oriented simple closed curve which is homologically nontrivial in $A(\Gamma)$.  The union of the sutures is denoted $s(\Gamma)$.

Every component of $R(\Gamma) = \partial Y - \mbox{Int}(\Gamma)$ is assigned an orientation compatible with the oriented sutures.  More precisely, if $\delta$ is a component of $\partial R(\Gamma)$, endowed with the boundary orientation, then $\delta$ must represent the same homology class in $H_1(\Gamma)$ as some suture.  Let $R_+(\Gamma)$ (resp., $R_-(\Gamma)$) denote those components of $R(\Gamma)$ whose normal vectors point out of (resp., into) $Y$.
\end{definition}

We omit mention of the distinguished annuli, $\Gamma$, from the notation when $Y$ is one of the following standard sutured manifolds:

\begin{example}\label{example:product} Let $F$ be an oriented surface with $\partial F \neq \emptyset$.  Then $F \times I$ denotes the {\em product sutured manifold} $(F \times I, \Gamma)$, where $\Gamma = \partial F \times I$, and $s(\Gamma) = \partial F \times \{\frac{1}{2}\}$.  %Note that $(F_{g,b},\boldalpha = \emptyset, \boldbeta = \emptyset)$ is a sutured Heegaard diagram for $F_{g,b} \times I$.
\end{example}

\begin{example}\label{example:branched} Let $(Y,\Gamma)$ be a sutured manifold and \[(B,\partial B) \subset (Y, \partial Y)\] a smoothly imbedded $1$--submanifold satisfying 
\[\partial B \cap \Gamma = \emptyset.\]
Let $\boldSigma(Y,B)$ be any two-fold cyclic branched cover of $Y$ over $B$ with covering map \[\pi: \boldSigma(Y,B) \rightarrow Y.\]

Then we denote by $(\boldSigma(Y,B), \widetilde{\Gamma})$ the sutured manifold with $\widetilde{\Gamma} = \pi^{-1}(\Gamma)$ and sutures $s(\widetilde{\Gamma}) = \pi^{-1}(s(\Gamma))$.  If $(Y,\Gamma)$ is a product sutured manifold, then we omit mention of $\widetilde{\Gamma}$ from the notation for $\boldSigma(Y,B)$.
\end{example}

\begin{remark}
Note that in general the construction of a $2$--fold cyclic branched cover of $Y$ over $B$ depends upon a choice of homomorphism \[\phi: H_1(Y-B;\Zint) \rightarrow \Ztwo.\]  %As long as $H_*(Y;\Zint)$ has no $2$--torsion (as is, for example, the case when $Y$ is a product sutured manifold $\clubsuit$), this is equivalent to a choice of cohomology class, $\phi \in H^1(Y-B;\Ztwo).$
In the present work, we make the following canonical choices:
\begin{enumerate}
  \item When $Y = D \times I$ and $B = T \subset D \times I$ is a balanced tangle with $m$ connected components, \[H_1(Y-B;\Zint) \cong \Zint^m,\] with basis given by the meridians $\mu_1, \ldots, \mu_m$.  We choose the homomorphism, $\phi$, satisfying $\phi(\mu_i) = 1 \in \Ztwo$ for all $\mu_i$.
\item When $Y = A \times I$ and $B = \bL \subset A \times I$ is a $0$--balanced tangle (link) with $m$ connected components, \[H_1(Y-B;\Zint) \cong \Zint^{m+1},\] with basis given by the meridians $\mu_1, \ldots, \mu_m$ along with $\gamma$, where $\gamma$ is represented by either connected component of $s(\Gamma)$.  We choose the homomorphism, $\phi$, satisfying $\phi(\mu_i) = 1$ for all $\mu_i$ and $\phi(\gamma) = 0$.  Note that, since the two components of $s(\Gamma)$ are homologous mod 2 in $Y-B$, the choice of $\gamma$ is irrelevant. 
\end{enumerate}
\end{remark}
%\begin{definition} \cite[Defn. 5.1-5.2]{GT08071432} \label{defn:balancedtangle} Let $F$ be an oriented surface with $\partial F \neq \emptyset$, and let $F \times I$ be the associated product sutured manifold, with $F_+:= F \times \{1\}$ (resp., $F_-:= F \times \{0\}$).  

%An admissible $n$--balanced tangle, $T \subset F \times I$ is (any representative of) an equivalence class of properly-imbedded (unoriented) $1$--manifolds satisfying

%\begin{itemize}
%  \item $T \cap \partial(F \times I) \subset \mbox{int} (F_+) \cup \mbox{int} (F_-)$,
%  \item $T_1, T_2$ are equivalent if they can be connected by an ambient isotopy acting trivially on $\partial F \times I$,
%  \item $\#(T \cap F_+) = \#(T \cap F_-) = n \in \Zint_{\geq 0}$.
%\end{itemize}

%We shall call any $0$--balanced tangle $\bL \subset F \times I$ a {\em link}, and any $1$--component link a {\em knot}.
%\end{definition}

\begin{definition}\cite[Defn. 2.7]{MR2390347} \label{defn:DecSurf} A {\em decomposing surface} in a sutured manifold, $(Y,\Gamma)$, is a properly imbedded, oriented surface, $(S,\partial S) \subset (Y,\partial Y)$ such that for every component, $\lambda$, of $(\partial S) \cap \Gamma$, one of the following occurs:
  \begin{itemize}
    \item $\lambda$ is a properly imbedded non-separating arc in $\Gamma$ such that $|\lambda \cap s(\Gamma)| = 1$.
    \item $\lambda$ is a simple closed curve in an annular component, $A$, of $\Gamma$ representing the same homology class in $A$ as $s(\Gamma)$.
    \item $\lambda$ is a homotopically nontrivial curve in a torus component $T$ of $\Gamma$, and if $\delta$ is another component of $T \cap (\partial S)$, then $\lambda$ and $\delta$ represent the same homology class in $T$.
  \end{itemize}
\end{definition}

\begin{definition} \label{defn:DecResult} Given a decomposing surface, $S$, in a sutured manifold, $Y$, the result of {\em decomposing $Y$ along $S$} is a new sutured manifold, $(Y',\Gamma')$, obtained as follows.
\begin{itemize}
  \item $Y' = Y - \text{Int}(N(S))$,
  \item $\Gamma' = (\Gamma \cap Y') \cup N(S_+' \cap R_-(\Gamma)) \cup N(S'_- \cap R_+(\Gamma))$,
  \item $R_+(\Gamma') = ((R_+(\Gamma) \cap Y') \cup S_+') - \text{Int}(\Gamma')$,
  \item $R_-(\Gamma') = ((R_-(\Gamma) \cap Y') \cup S_-') - \text{Int}(\Gamma')$,
\end{itemize}
where $S_+'$ (resp., $S'_-$) is the component of $\partial N(S) \cap Y'$ whose normal vector field points out of (resp., into) $Y'$. 
\end{definition}

\begin{definition}
A {\it link} $\bL \subset A \times I$ is a smoothly imbedded, unoriented, closed $1$--manifold. Two such imbeddings $\bL_1, \bL_2$ are said to be equivalent if there is an ambient isotopy taking $\bL_1$ to $\bL_2$ which acts trivially on $\partial A \times I$.  A {\it knot} is a $1$--component link.
\end{definition}

Note that we require our ambient isotopy to act trivially on $\Gamma = \partial A \times I$, the neighborhood of the sutures, in order to preserve the structure of the ambient manifold, $A \times I$ as a sutured manifold.
%\begin{definition}
%A {\it framed knot} $(K,\lambda) \subset A \times I$ is a knot $K \subset A \times I$ along with a smooth section, $\lambda$, of its normal disk bundle, $N_{A \times I}(K)$, satisfying $\lambda \cap K = \emptyset$.  Such a pair determines a smoothly imbedded annulus $A_{K,\lambda} \subset A \times I$.  A {\it framed link} $(\bL, \vec{\lambda}) = (\vec{K}, \vec{\lambda})$ is (the equivalence class of) a disjoint union, \[A_{\bL,\vec{\lambda}} = \coprod_{i=1}^n A_{K_i, \lambda_i},\] of such imbeddings, where two framed links, $(\bL_1,\vec{\lambda}_1)$ and $(\bL_2, \vec{\lambda}_2)$ are said to be equivalent if there is an ambient isotopy connecting $A_{\bL_1,\vec{\lambda}_1}$ to $A_{\bL_2,\vec{\lambda}_2}$ which acts trivially on $\partial A \times I$.
%\end{definition}

\begin{definition}
Let $\pi_z: A \times I \rightarrow A$ be the standard projection defined by \[\pi_z(r,\theta,z) = (r,\theta).\]  For any link $\bL \subset A \times I$ for which $\pi_z(\bL) \subset A$ is a smooth imbedding away from finitely many transverse double points, we denote by $\cP(\bL)$ the enhancement of $\pi_z(\bL)$ which encodes over/undercrossing information.  We call $\cP(\bL)$ the {\em projection} of $\bL$.
\end{definition}

%\begin{remark} \label{rmk:framing} A projection, $\cP(\bL)$, specifies a {\em blackboard framing} %$\vec{\lambda}$, 
%for $\bL$ in the standard way.  %as follows.  $\pi_z(\bL)$ is an immersed multi-curve in $A$, hence has a normal bundle, $N_A(\cP(\bL))$.  The blackboard framing of $\bL$ induced by $\cP(\bL)$ is a smooth section, $\vec{\lambda}$, of the normal bundle of $\bL \subset A \times I$ satisfying the property that $\pi_z(\lambda)$ is a smooth section of $N_A(\cP(\bL))$ which intersects $\cP(\bL)$ exactly once in the neighborhood of each transverse double point, each time in the complement of the double point.  In this way
%Therefore, a link projection $\cP(\bL)$ uniquely determines not only a link equivalence class but also a framed link equivalence class.
%\end{remark}

%\begin{definition} A (framed) link equivalence class representative $\bL$ is said to be {\em resolved} if $\cP(\bL)$ is an imbedding. \end{definition}

\begin{definition} A link $\bL \subset A\times I$ is said to be {\em resolved} if $\cP(\bL)$ is an imbedding. \end{definition}

\begin{definition}
A {\em saddle cobordism} $S\subset A\times [0,1]$ is a smooth cobordism between two resolved link projections $\cP(\bL')$ and $\cP(\bL'')$ with the property that $\exists$ a unique $c \in [0,1]$ such that 

\begin{enumerate}
\item $S \cap (A \times \{c\})$ is a smooth $1$--dimensional imbedding away from a single double-point.
\item $S \cap (A \times \{s\})$ is a smooth $1$--dimensional imbedding whenever $s \neq c$.
\end{enumerate}

Let $|\bL'|$ (resp., $|\bL''|$) denote the number of connected components of $\bL'$ (resp., $\bL''$).  There are two cases:

\begin{enumerate}
  \item When $|\bL'| = |\bL''| + 1$, we call $S$ a {\em merge saddle cobordism}, and
  \item when $|\bL'| = |\bL''|-1$, we call $S$ a {\em split saddle cobordism}.
\end{enumerate}
\end{definition}

\subsection{Khovanov functor} \label{sec:KBSM}
%For a resolved link $\bL\subset A\times I$ with $t$ trivial
%components $K_1,\ldots,K_t$ and $n$ non-trivial components
%$K_{t+1},\ldots,K_n$, 
%let $Z(\bL)$ be the $\Ztwo$--vector space formally generated
%by all components $[K_1],\ldots,[K_{t+n}]$, and let
%$V(\bL)$ denote the exterior algebra,
%$$V(\bL):=\Wedge^*Z(\bL),$$ i.e. this
%is the polynomial algebra over $\mathbb{Z}_2$ in formal variables
%$[K_1],\ldots,[K_{t+n}]$, modulo the relations $[K_i]^2=0$
%for $i=1,\ldots,t+n$.

Let $\bL \subset A \times I$ be a %(framed)
resolved link with connected components $K_1, \ldots, K_t, K_{t+1}, \ldots, K_{t+n}$, where \[\cP(K_i) \,\,\,\, \left\{\begin{array}{ll} =0\in H_1(A;\Zint)  & \mbox{ if $i \in \{1, \ldots, t\},$ and}\\
          \neq 0\in H_1(A;\Zint) & \mbox{ if $i \in \{t+1, \ldots, t+n\}.$}\end{array}\right.\]  Denote 
\begin{eqnarray*}
\bL_{triv} &:=& K_1 \amalg \ldots \amalg K_t\\
\bL_{non} &:=& K_{t+1} \amalg \ldots \amalg K_{t+n}.
\end{eqnarray*}

%It will also be convenient for us to define/denote:
%\begin{itemize}
%  \item $D \times I := \{(r,\theta,z) \in \Real^3\,\,|\,\, r \in [0,2],z \in [0,1]\}$,
%  \item $i: A \times I \longrightarrow D \times I$ the inclusion map,
%  \item $\hatbL := i(\bL) \subset D \times I$ the image of $\bL$ under the inclusion.
%\end{itemize} 

\begin{definition}
For $\bL = K_1 \amalg \ldots \amalg K_{t+n}\subset A \times I$ a resolved link, let $Z({\bL})$ denote the $\Ztwo$--vector space formally generated by $[K_1], \ldots, [K_{t+n}]$: \[Z(\bL) := Span_{\Ztwo}([K_1], \ldots, [K_{t+n}]).\]
Define a bigrading on $Z(\bL)$ by
$$
deg([K_i]):=
\begin{cases}
(0,-2) & \text{for $1\leq i\leq t$, and}\\
(-2,-2)& \text{else.}
\end{cases}
$$
The space $Z(\bL)$ decomposes as
$Z(\bL)=Z(\bL_{triv})\oplus Z(\bL_{non})$, where
\begin{eqnarray*}
Z(\bL_{triv}) &:=& Span_{\Ztwo}([K_1], \ldots, [K_{t}]).\\
Z(\bL_{non}) &:=& Span_{\Ztwo}([K_{t+1}], \ldots, [K_{t+n}]).
\end{eqnarray*}
\end{definition}

\begin{definition}\label{defn:VofL}
Let $\overline{V}(\bL)$ denote the exterior algebra
$$
\overline{V}(\bL) := \Wedge^*Z(\bL),
$$
i.e. the polynomial algebra over $\mathbb{Z}_2$ in
formal variables $[K_1], \ldots, [K_{t+n}]$ satisfying the relations
$[K_1]^2\sim\ldots\sim[K_{t+n}]^2\sim0$.
The bigrading on $Z(\bL)$ induces a bigrading
on $\overline{V}(\bL)$, via
$$
deg([K_{i_1}]\wedge\ldots\wedge[K_{i_k}]):=
deg([K_{i_1}])+\ldots+deg([K_{i_k}]).
$$
Given a bigraded vector space $V$ and a pair $(a,b)\in\mathbb{Z}^2$,
we denote by $V\{a,b\}$ %\footnote{$\clubsuit$ Check that this is the standard convention for grading shifted modules.} 
the vector space $V$, with
gradings shifted by $(a,b)$.  I.e., $deg_{V\{a,b\}}(v):=deg_V(v) + (a,b)$
for all $v\in V$.
Let
$$
V(\bL):=\overline{V}(\bL)\{n,t+n\}.
$$
We can write $V(\bL)$ as a direct sum
$$
V(\bL)=\bigoplus_{f,q\in\mathbb{Z}}V(\bL;f,q)
$$
where $V(\bL;f,q)$ denotes the subspace
consisting of all $v\in V(\bL)$ with $deg(v)=(f,q)$.
\end{definition}
%The vector space $V(\bL)$ comes with a bigrading:
%$$
%V(\bL)=\bigoplus_{f,q\in\mathbb{Z}} V(\bL,f,q)
%$$
%where $V(\bL,f,q)\subset\Wedge^*Z(\bL)=V(\bL)$ denotes the subspace 
%$$ V(\bL,f,q):=\Wedge^{\frac{t + (q-f)}{2}} Z(\bL_{triv})
%\wedge\Wedge^{\frac{n+f}{2}}
%Z(\bL_{non}),
%\footnote{$\clubsuit$ I changed this from what you had originally
%written, Stephan.  I think this is correct, but you should check it.
%What you had originally written was:
%$V(\bL,f,q):=\Wedge^{t-2(q-f)} Z(\bL_{triv})
%\wedge\Wedge^{n-2f} Z(\bL_{non})$.}
%$$
%subject to the rule  that $\Wedge^a(Z) = 0$ if $a \not\in \Zint$.
%\end{definition}

Now consider a merge saddle cobordism $S_m\subset A\times [0,1]$
between two resolved link projections $\cP(\bL')$ and $\cP(\bL'')$,
where the saddle merges two components of $\bL'$ labeled $K'_i$ and $K'_j$.
Then there is a natural identification
$$
Z(\bL')/[K'_i]\sim [K'_j]=Z(\bL'')
$$
and, correspondingly, there is an isomorphism
$$
\alpha : V(\bL')/[K'_i]\sim [K'_j]\stackrel{\cong}{\longrightarrow} V(\bL'').
$$

\begin{definition}\label{defn:KhMergeMap}
Associated to $S_m$ is a linear
map $\cV_m:V(\bL') \rightarrow V(\bL'')$,
referred to as the multiplication map, and
defined as the composite
\[
\begin{CD} V(\bL') @>{\pi}>> 
\frac{V(\bL')}{[K'_i]\sim [K'_j]
}@>{\alpha}>> V(\bL''),
\end{CD}
\]
where $\pi$ denotes the quotient map.  Let $G\cV_m$ denote the part of $\cV_m$ which is $f$--grading preserving on $V(\bL')$, $V(\bL'')$.  More explicitly, if ${\bf x} \in V(\bL'')$ is a homogeneous element in $f$--grading $f_0$, then \[G\cV_m({\bf x}) := p_{f_0}\circ\cV_m({\bf x}),\]  where $p_{f_0}$ is the projection map $p_{f_0}: V(\bL'') \rightarrow V(\bL'';f=f_0)$.  Now extend linearly.
\end{definition}

Running the merge saddle cobordism $S_m$ backwards
produces a split saddle cobordism
$S_{\Delta}:\cP(\bL'')\rightarrow\cP(\bL')$. 

\begin{definition} \label{defn:KhSplitMap}
Associated
to $S_{\Delta}$ is a comultiplication map
$\cV_{\Delta}:V(\bL'') \rightarrow V(\bL')$,
defined as the composite
$$
\begin{CD}V(\bL'') @>{\alpha^{-1}}>> 
\frac{V(\bL')}{[K'_i]\sim [K'_j]}@>{\varphi}>> V(\bL'),
\end{CD}
$$
with $\varphi$ given by
$\varphi(a):=([K'_i]+[K'_j])\wedge \widetilde{a}$, where $\widetilde{a}$ is any lift of $a$ in $\pi^{-1}(a)$.  As above, we denote by $G\cV_\Delta$ the part of $\cV_\Delta$ which is $f$--grading perserving, i.e., $G\cV_\Delta$ linearly extends the map which sends elements ${\bf x} \in V(\bL'')$ in $f$--grading $f_0$ to \[p_{f_0} \circ \cV_\Delta({\bf x}).\]
\end{definition}

The maps %$\cV_m$, $\cV_{\Delta}$ and 
$G\cV_m$, $G\cV_\Delta$ are used to define chain complexes associated
to a projection $\cP(\bL)$ of a link
$\bL\subset A\times I$ as follows.
   
Label the crossings of $\cP(\bL)$ by $1, \ldots, \ell$. 
For any $\ell$--tuple $\cI = (m_1, \ldots, m_{\ell}) \in \{0,1,\infty\}^\ell$,
we denote by $\cP_{\cI}(\bL)$ the link projection obtained from $\cP(\bL)$
by

\begin{itemize}
  \item leaving a neighborhood of the $i$th crossing unchanged, if $m_i = \infty$,
  \item replacing a neighborhood of the $i$th crossing with a ``0'' resolution, if $m_i = 0$, and
  \item replacing a neighborhood of the $i$th crossing with a ``1'' resolution, if $m_i = 1$.
\end{itemize}

\begin{figure}
\begin{center}
\resizebox{2in}{!}{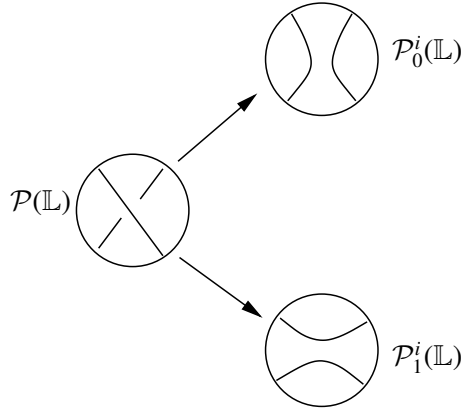}
\end{center}
\caption{Replacing the $i$-th crossing of a projection, $\cP(\bL)$, with a $0$ or $1$ resolution}
\label{fig:Resolutions}
\end{figure}

See Figure \ref{fig:Resolutions}.  Note that our conventions for the ``0'' and ``1'' resolutions of a crossing match those of \cite{MR2141852} and \cite{GT08071432}, which are opposite of Khovanov's standard convention \cite{MR1740682}, used by Roberts in \cite{GT07060741}.

%It is clear that
%$\cP_{\cI}(\bL)$ is resolved (no crossings)
%whenever $\cI$ is contained in $\{0,1\}^{\ell}$. Note also that
%$\cP_{\cI}(\bL)$ and $\cP_{\cI'}(\bL)$
%related to each other by a saddle cobordism
%whenever $\cI'\in\{0,1\}^{\ell}$ is an immediate successor
%of $\cI\in\{0,1\}^{\ell}$. We define:

\begin{definition} Giving the set $\{0,1,\infty\}$ the dictionary ordering, we call an $\ell$--tuple $\mathcal{I'} \in \{0,1,\infty\}^\ell = (m_1, \ldots, m_\ell)$ an {\em immediate successor} of $\mathcal{I}$ if there exists some $j$ such that $m_i = m_i'$ if $i \neq j$ and $(m_j, m_j')$ is either $(0,1)$ or $(1,\infty)$.  
\end{definition}

\begin{definition}\label{defn:KhChainCpx}
Given a projection $\cP(\bL) \subset A$ of a link
$\bL\subset A\times I$,
we define a chain complex \[CV(\cP(\bL)) = \left(\bigoplus_{\cI \in \{0,1\}^\ell} V(P_{\cI}(\bL)),D\right),\] where $D = \sum_{\cI,\cI'} D_{\cI,\cI'}$, with the sum taken over all pairs $\cI,\cI' \in \{0,1\}^\ell$ such that $\cI'$ is an immediate successor of $\cI$, and \[D_{\cI, \cI'}: V(\cP_{\cI}(\bL)) \rightarrow V(\cP_{\cI'}(\bL))\] is
given by $\cV_m$ (resp., $\cV_\Delta$) when $\cP_{\cI'}(\bL)$ is obtained from $\cP_{\cI}(\bL)$
by a merge (resp., split) saddle cobordism.
\end{definition}

If $\bL$ is equipped with an orientation,
then the complex $CV(\cP(\bL))$
can be endowed with three $\mathbb{Z}$ gradings,
called the $i$, $j$ and $k$--grading.

\begin{definition}\label{defn:KhGradings}
An element $v\in CV(\cP(\bL))$ is said to have
degree $deg(v)=(i,j,k)\in\mathbb{Z}^3$ if
it is contained in the subspace
$$
CV(\cP(\bL);i,j,k):=\bigoplus V_{\cI}(\cP(\bL);f,q)
\;\subset\; CV(\cP(\bL))
$$
where the sum ranges over all
triples $(\cI,q,f)\in\{0,1\}^{\ell}\times\mathbb{Z}^2$ satisfying
\begin{align*}
i& = |\cI|-n_+,\\
j&= q + |\cI|+n_- -2n_+,\\
k&= f
\end{align*}
where $n_+$ (resp., $n_-$) denotes the
number of positive (resp., negative) crossings
in $\cP(\bL)$, and $|\cI|$ denotes the number
of $1$'s in the $\ell$--tuple $\cI\in\{0,1\}^{\ell}$.
\end{definition}

%If $\bL$ is equipped with an orientation,
%then the complex $CV(\cP(\bL))$ can be endowed
%with three $\mathbb{Z}$
%gradings, as follows.
%%called the $i$, $j$ and $k$--grading.
%%Let $n_+$ (resp., $n_-$) denote the number of positive
%%(resp., negative) crossingss in $\bL$).
%For an
%$\ell$--tuple $\cI=(m_1,\ldots,m_{\ell})\in\{0,1\}^{\ell}$, let
%$|\cI|:=\sum_i m_i$
%$$
%CV(\cP(\bL),i):=
%\left(\bigoplus_{|\cI|-n_-=i}V(\cP_{\cI}(\bL))\right)\{n_--2n_++i,0\}
%$$
%where $n_+$ (resp., $n_-$) is the number of positive (resp., negative)
%crossings in $\cP(\bL))$, and $\{n_--2n_++i,0\}$ is the functor
%which raises the $(f,q)$ bigrading of the $V(\cP_{\cI}(\bL))$ by $(n_--2n_++i,0)$,
%as in Definition~\ref{defn:VofL}.

%{\it $\clubsuit$ Stephan, we should discuss the i,j,k gradings at some point, so that we can verify that my ideas about them match yours.}

%\begin{definition}({\em $i$, $j$ and $k$--gradings})\label{defn:KhGradings}
%Let ${\bf x}$ be a generator of $CV(\cP(\bL))$
%satisfying ${\bf x} \in V_{\cI}(\cP(\bL),q,f)$,
%for an $\ell$--tuple
%$\cI=(m_1,\ldots,m_{\ell}) \in\{0,1\}^{\ell}$,
%and integers $q,f\in\mathbb{Z}$. Then
%\begin{align*}
%i({\bf x})&:= -n_++\sum_i m_i,\\
%j({\bf x})&:= q + n_- -2n_+ + \sum_i m_i,\\
%k({\bf x})&:= f.
%\end{align*}
%where $n_+$ (resp. $n_-$) denotes the number of
%positive (resp. negative) crossings in $\cP(\bL)$.
%\end{definition}

\begin{remark} \label{rmk:fgrading} For a resolved link, $\bL \subset A \times I$, our $V(\bL)$ is equivalent, as a $\Ztwo$ vector space, to what Roberts denotes $V(\bL)$ in
\cite[Sec. 2]{GT07060741}, and our $(f,q)$--bigrading matches his for resolved links.  For an unresolved link $\bL$, Roberts' $f$--grading (for $\overline{\bL}$) matches ours (for $\bL$), and his $q$--grading matches our $j$--grading (which in turn matches our $q$--grading up to an overall shift). Note also that, if one forgets the $k$--grading, then $CV(\cP(\bL))$ becomes the chain complex that Ozsv\'ath and Szab\'o associate to $\bL \subset S^3$ in \cite{MR2141852}. In particular, the homology of $CV(\cP(\bL))$ is the Khovanov homology of $\overline{\bL} \subset S^3$.

%\footnote{$\clubsuit$ There is still the issue that for unresolved links, Roberts' $q$ grading for the mirror agrees with our $q$ grading for the original link.  We need to address this.  We also need a remark after the definitions of $0$, $1$ resolutions on the next page.  Perhaps move this remark later and address both at the same time?}   Compare Definitions \ref{defn:KhChainCpx} and \ref{defn:KhGradings} below to \cite{GT07060741}.
\end{remark}

\begin{lemma}
The differential, $D$, on $CV(\cP(\bL))$, is non-increasing in the $k$ grading.  Hence, $CV(\cP(\bL))$ has the structure of a $\Zint$--filtered chain complex, with the $\Zint$--filtration induced by the $k$--grading.
\end{lemma}

\begin{proof} Let $\cI, \cI' \in \{0,1\}^\ell$, where $\cI'$ is an immediate successor of $\cI$, and let $\cP_{\cI}(\bL), \cP_{\cI'}(\bL)$ denote the associated resolutions.  

Suppose that $\cP_{\cI'}(\bL)$ is obtained from $\cP_\cI(\bL)$ by merging two components, $K_i$ and $K_j$, of $\cP_\cI(\bL)$ to form a component, $K'$, of $\cP_{\cI'}(\bL)$.  Let $deg_f(K_i), deg_f(K_j)$ (resp., $deg_f(K')$) denote their $f$--degrees as elements of $Z(\cP_\cI(\bL))$ (resp., $Z(\cP_{\cI'}(\bL))$).

%Note that \[deg([K_i']) = \min\{deg([K_i]),deg([K_j])\},\] where \[\min\{(f_1,q_1),(f_2,q_2)\} := \left(\min\{f_1,f_2\},\min\{q_1,q_2\}\right)\]

There are three possibilities:
\begin{enumerate}
  \item $K_i, K_j \subset \left(\cP_\cI(\bL)\right)_{triv}$% merge to form a component, $K_i' \subset \bL_{triv}(\cI')$,
  \item one of $K_i, K_j$ is in $\left(\cP_\cI(\bL)\right)_{triv}$, while the other is in $\left(\cP_\cI(\bL)\right)_{non}$, or % merges with a component of $\bL_{triv}(\cI)$ to form a component of $\bL_{non}$, or 
  \item $K_i, K_j \subset \left(\cP_\cI(\bL)\right)_{non}$.% merge to form a component of $\bL_{triv}(\cI')$.
\end{enumerate}

%In cases $1$ and $3$ above, $K_i$ and $K_j$ merge to form a component of $\bL_{triv}(\cI')$, while in case $2$, they merge to form a component of $\bL_{non}(\cI')$.  

In the first two cases, \[deg_f(K') = \min\left\{\,deg_f(K_i)\, ,\,deg_f(K_j)\,\right\},\]  and since, in these two cases, $\cP_\cI(\bL)$ and $\cP_{\cI'}(\bL)$ have the same number of non-trivial components, the $f$--grading shift relating $V\left(\cP_\cI(\bL)\right)$ to  $\overline{V}\left(\cP_\cI(\bL)\right)$ is the same as the $f$--grading shift relating $V\left(\cP_{\cI'}(\bL)\right)$ and $\overline{V}\left(\cP_{\cI'}(\bL)\right)$.  Hence, $D_{\cI,\cI'}:=\alpha \circ \pi$ will be non-increasing (specifically, either degree $0$ or $-2$) in the $(k=f)$--grading.

In the third case, $deg_f(K') = 0$, while $deg_f(K_i) = deg_f(K_j) = -2$.  Since $\cP_{\cI'}(\bL)$ has two fewer non-trivial components than $\cP_\cI(\bL)$, after applying the $f$--grading shift to $\overline{V}\left(\cP_\cI(\bL)\right)$ (resp., $\overline{V}\left(\cP_{\cI'}(\bL)\right)$) to obtain $V\left(\cP_\cI(\bL)\right)$ (resp., $V\left(\cP_{\cI'}(\bL)\right)$) as above, we similarly conclude that $D_{\cI,\cI'}$ is non-increasing in the $(k=f)$--grading.

%Note that, considered as elements of $Z(\bL)$, $Z(\bL')$ in cases $1$ and $2$, we have the grading on $(\bL(\cI'))$ is induced from $V(\bL(\cI))$ map $D_{\cI, \cI'} = \alpha \circ \pi$ and, hence, is automatically non-increasing in the $k$--grading.  %This is because the $k$ grading of any generator of $V(\bL(\cI'))$ is the sum of the $k$ gradings, and the merge map is like a multiplication map.
%In the third case, the number of components of $\bL_{non}(\cI')$ is $2$ fewer than the number of components of $\bL_{non}(\cI)$, again ensuring that $D_{\cI,\cI'} = \alpha \circ \pi$ either preserves the $k$--grading or decreases it by $2$.  $\bL_{non}(\cI)$$\bL_{non}(\cI')$ In the first case, $D_{\cI,\cI'}:= \alpha \circ \pi$ preserves the $k$ grading on all generators.  In the second case, $\alpha \circ \pi$ is 

%If , then $D_{\cI,\cI'}$ is $k$--grading preserving.  and let ${\bf x} \in V(\cP_\cI(\bL);f_x, q_x)$.  Suppose ${\bf y} \in V(\cP_{\cI'}(\bL);f_y,q_y)$ is in the support of $D_{\cI,\cI'}$.  Then \[q_y - q_x = \left\{\begin{array}{rl}
%                     0 & \mbox{if $D_{\cI,\cI'} = \cV_m$, and}\\
%		     -2 & \mbox{if $D_{\cI,\cI'} = \cV_\Delta$.}
%		   \end{array}.\]
The case where $\cI'$ is obtained from $\cI$ by a split saddle cobordism is completely analogous.  
\end{proof}

\begin{definition}\label{defn:GrKhChainCpx}
Let $GCV(\cP(\bL))$ denote the associated graded complex of the $k$--filtered complex $CV(\cP(\bL))$.  I.e., 
\[GCV(\cP(\bL)) \cong \left( CV(\cP(\bL)) ,GD\right),\] where $GD= \sum_{\cI,\cI'} GD_{\cI,\cI'}$, with the sum taken over all pairs $\cI,\cI' \in \{0,1\}^\ell$ where $\cI'$ is an immediate successor of $\cI$ and \[GD_{\cI,\cI'}: V(\cP_{\cI}(\bL)) \rightarrow V(\cP_{\cI'}(\bL))\] is given by $G\cV_m$ (resp., $G\cV_\Delta$) when $\cP_{\cI'}(\bL)$ is obtained from $\cP_{\cI}(\bL)$ by a merge (resp., split) saddle cobordism. Let $V(\bL)$ denote the homology of $GCV(\cP(\bL))$.   
\end{definition}

\begin{remark}\label{rmk:FrameDep}
Our chain complex $GCV(\cP(\bL))$ is isomorphic to the triply--graded chain complex $C^{*;*,*}(\cP(\overline{\bL}))$ that Roberts associates to $\bL \subset A\times I$ in \cite{GT07060741}. Roberts' complex, in turn, is essentially equivalent to the Asaeda-Przytycki-Sikora chain complex for $\cP(\bL)$, $\bL \subset A \times I$, defined in \cite{MR2113902}, except that the Asaeda-Przytycki-Sikora construction requires $\bL$ to be framed (but not oriented), whereas Roberts' construction requires $\bL$ to be oriented (but not framed). Note that Asaeda-Przytycki-Sikora's $s$--grading appears (without loss of information) as Roberts' $f$--grading and that the relative $i$-- and $j$--grading on $GCV(\cP(\bL))$, as well as the absolute $k$--grading, do not depend on the choice of an orientation for $\bL$.  Finally, we remark that the Khovanov-type homology of $\bL$, $Kh^*(\bL)$, referenced in the introduction and in the statements of Theorems \ref{thm:Reinterp} and \ref{thm:Cutting}, is $V(\overline{\bL})$.  See Remark \ref{rmk:fgrading}.

\end{remark}

%Subject to Remark \ref{rmk:FrameDep}, we make the following definitions.

%\begin{definition}\label{defn:FramIndepDef}
%Let $\bL \subset A \times I \subset S^3$ be a link and let $\cP^*(\bL)$ be any projection of $\bL$ inducing a blackboard framing which matches the Seifert framing on each component.  

%Then we denote by $CV(\bL)$ the chain complex, $CV(\cP^*(\bL))$, of Definition \ref{defn:KhChainCpx} and by $GCV(\bL)$ the chain complex, $GCV(\cP^*(\bL))$, of Definition \ref{defn:GrKhChainCpx}.
%\end{definition}

%\begin{definition} (Quantum grading)
%Let $a$ be an element
%of $CV(\cP(T))$,
%and suppose that $a$ is contained in
%$\Wedge^d (Z(\cP_{\cI}(T)))\subset V(\cP_{\cI}(T))
%\subset CV(\cP{T})$, for an
%$\ell$-tuple $\cI=(m_1,\ldots,m_\ell)\in\{0,1\}^\ell$
%and a non-negative integer $d\geq 0$.
%Then
%where $n_+$ (resp. $n_-$) denotes the
%number of positive (resp. negative) crossings
%in $\cP(T)$.
%\end{definition}

%{\it $\clubsuit$ We may need a few words to segue to the next subsection here.}

\subsection{Sutured Floer Functor} \label{sec:SFH}
Given a link $\bL \subset A \times I$, we alternatively obtain a filtered chain complex using sutured Floer homology (see \cite{MR2253454}), a theory developed by Juh{\'a}sz which assigns Floer homology groups to balanced sutured manifolds.  Recall (\cite[Defn. 2.2]{MR2253454}) that a sutured manifold $(Y,\Gamma)$ is said to be {\em balanced} if $\chi(R_+) = \chi(R_-)$, and the maps $\pi_0(\Gamma) \rightarrow \pi_0(\partial Y)$ and $\pi_0(\partial Y) \rightarrow \pi_0(Y)$ are surjective. Given a balanced sutured manifold $(Y,\Gamma)$--for example the branched double-cover, $\boldSigma(A\times I,\bL)$, where $\bL\subset A\times I$ is a link as above--one obtains its sutured Floer homology, $SFH(Y,\Gamma)$, by the following procedure:

%Since $\bL \subset A \times I$ satisfies $\#(\bL \cap A_+) = \#(\bL \cap A_-)$, $\boldSigma(A \times I, \bL)$ is balanced (\cite[Prop. 2.12]{GT08071432}).  Hence its {\em sutured Floer homology}, $SFH(\boldSigma(A \times I, \bL))$, can be defined as in \cite{MR2253454}.  We summarize the procedure below for the convenience of the reader.
%  Before describing the sutured Floer functor, we will need to summarize some basic defintions and results.
%\begin{definition} \cite[Defn. 2.2]{MR2253454} \label{defn:balancedSM} 
%\end{definition}

%\begin{enumerate}
%  \item $A \times I$ is naturally a {\em balanced, sutured manifold} (\cite[Defn. 2.2]{MR2253454}) with $\Gamma:= (\partial A) \times I$ and $s(\Gamma) = \partial A_{\frac{1}{2}}$, and
%  \item $\bL \subset A \times I$ satisfies $\#(\bL \cap A_+) = \#(\bL \cap A_-)$; hence, $\boldSigma(A \times I, \bL)$ is also a balanced, sutured manifold (\cite[Prop. 2.12]{GT08071432}).
%\end{enumerate}

%We briefly recall the relevant definitions here, referring the interested reader to \cite{MR2253454} for more details.

\begin{enumerate}
     \item Construct a {\em balanced, sutured Heegaard diagram} (\cite[Defn. 2.7, 2.8]{MR2253454}), $(\Sigma, \boldalpha, \boldbeta)$, for $(Y,\Gamma)$, where $\Sigma$ is a compact, oriented surface with non-empty boundary and no closed components, and $\boldalpha = \{\alpha_1, \ldots, \alpha_d\}$, $\boldbeta = \{\beta_1, \ldots, \beta_d\}$ are two sets of pairwise disjoint simple closed curves in Int$(\Sigma)$ satisfying the condition that the set $\boldalpha$ (resp., $\boldbeta$) is linearly independent in $H_1(\Sigma;\Zint)$.  Recall that one obtains a unique balanced, sutured manifold from a balanced, sutured Heegaard diagram by attaching $3$--dimensional $2$--handles to $\Sigma \times I$ along the curves $\alpha_i \times \{0\}$ and $\beta_j \times \{1\}$ for $i,j \in \{1, \ldots d\}$.  $\Gamma$ is $\partial \Sigma \times I$, and $s(\Gamma) = \partial \Sigma \times \{\frac{1}{2}\}$.

%\begin{definition} \cite[Defn. 2.7, 2.8]{MR2253454} \label{defn:suturedHD} A {\em sutured Heegaard diagram} is a tuple $(\Sigma,\boldalpha,\boldbeta)$, where $\Sigma$ is a compact, oriented surface with boundary, and $\boldalpha = \{\alpha_1, \ldots, \alpha_d\}$, $\boldbeta = \{\beta_1, \ldots, \beta_d\}$ are two sets of pairwise disjoint simple closed curves in Int$(\Sigma)$.  Every sutured Heegaard diagram uniquely defines the sutured manifold obtained by attaching $3$--dimensional $2$--handles to $\Sigma \times I$ along the curves $\alpha_i \times \{-1\}$ and $\beta_j \times \{4\}$ for $i,j \in \{1, \ldots d\}$.  $\Gamma$ is $\partial \Sigma \times I$, and $s(\Gamma) = \partial \Sigma \times \{\frac{3}{2}\}$.
%\end{definition}

%\Juhasz then defines sutured Floer homology for {\em balanced} sutured manifolds (see \cite[Defn. 2.2]{MR2253454}), which can be described via {\em balanced} Heegaard diagrams:

%\begin{definition} \cite[Defn. 2.1]{MR2253454} \label{defn:balancedHD}A sutured Heegaard diagram $(\Sigma, \boldalpha, \boldbeta)$ is called {\em balanced} if $|\boldalpha| = |\boldbeta|$, $\Sigma$ has no closed components, and $\{\alpha_i\}$ (resp., $\{\beta_i\}$) are linearly-independent in $H_1(\Sigma)$.
%\end{definition}

\item The data of a balanced Heegaard diagram \[(\Sigma, \boldalpha = \{\alpha_1, \ldots, \alpha_d\}, \boldbeta = \{\beta_1, \ldots, \beta_d\})\] and a generic (family of) complex structures on $\Sigma$, is used to construct a Floer chain complex using the half-dimensional tori $\Torus_\alpha = \alpha_1 \times \ldots \times \alpha_d$ and $\Torus_\beta = \beta_1 \times \ldots \times \beta_d$ in $Sym^d(\Sigma)$.  Specifically, one obtains a chain complex with:

\begin{enumerate}
\item Generators: $\{{\bf x} \in \Torus_\alpha \cap \Torus_\beta\}$,
\item Differentials: \[\partial({\bf x}) = \sum_{{\bf y} \in \Torus_\alpha \cap \Torus_\beta} \sum_{\{\phi \in \pi_2({\bf x},{\bf y})| \mu(\phi) = 1\}} \widehat{\cM}(\phi) \cdot {\bf y}.\]
\end{enumerate}

As usual, $\pi_2({\bf x},{\bf y})$ denotes the homotopy classes of disks (disjoint from $\partial \Sigma$) connecting ${\bf x}$ to ${\bf y}$, $\mu(\phi)$ denotes the Maslov index of $\phi \in \pi_2({\bf x},{\bf y})$, and $\widehat{\cM}(\phi)$ denotes the moduli space of holomorphic representatives of $\phi$, modulo the standard $\Real$ action.  

\item Denote by $CFH(Y,\Gamma)$ the chain complex constructed in (2), and by $SFH(Y,\Gamma)$ its homology.

%Denote by $CFH(Y,\Gamma)$ any chain complex associated to a balanced, sutured manifold $(Y,\Gamma)$ arising as above, and by $SFH(Y,\Gamma)$ the homology of such a chain complex.

\end{enumerate}

In preparation for describing a sutured Floer-type ``functor'' analogous to the Khovanov-type ``functor'' defined in Section \ref{sec:KBSM}, we now establish, in Proposition \ref{prop:SFHFiltCpx}, %and \ref{prop:SFHFiltCpx}, 
a relationship between the sutured Floer homology of $\boldSigma(A \times I, \bL)$ and the knot Floer homology of a certain link in $\boldSigma(S^3,\bL)$.  We shall see that $SFH(\boldSigma(A \times I, \bL))$ can be identified as the homology of the associated graded complex of a natural filtered complex for this link in $\boldSigma(A \times I,\bL)$.   Proposition \ref{prop:SFHFiltCpx} provides a sutured-Floer theoretic interpretation of the Floer homology invariants studied in \cite{GT07060741}.  

%After establishing this relationship, we will be ready to describe a sutured Floer-type ``functor'' analogous to the Khovanov-type ``functor'' defined in Section \ref{sec:KBSM}.

In what follows:
\begin{itemize}
%\item $S^3 := (\Real^2 \times \Real) \,\,\cup\,\, \infty$ and 
\item $A \times I = \{(r,\theta, z) \,|\, r \in [1,2], \theta \in [0,2\pi), z \in [0,1]\} \subset \Real^3 \cup \infty = S^3$ denotes the standard imbedding of $A \times I$ into $S^3$, with $\bL \subset A \times I$ a fixed (isotopy class of) link,
\item $B = \{(r,\theta, z)\,|\,r=0\} \cup \infty \subset S^3$ denotes a standardly-imbedded unknot,% isotopic in the complement of $\bL \subset S^3$ to $\partial(\lambda \times I)$
\item $\pi:\boldSigma(S^3,\bL) \rightarrow S^3$ denotes the branched-covering projection,
%\item $B := \{\mbox{z-axis} \,\cup \,\infty\}\sim \partial(D_0),$
\item  $\widetilde{B} := \pi^{-1}(B)$ denotes the preimage of $B$ in $\boldSigma(S^3, \bL)$,
%\item $\lambda := \{(r,\theta)\,|\,\theta = 0\} \subset A$ be a properly-imbedded radial arc in $A$, oriented in the direction of decreasing $r$,
%\item $\lambda \times I$ be the associated vertical disk, with the induced product orientation, and
  \item $D_A \times I := \{(r,\theta,z) \in \Real^3\,\,|\,\, r \in [0,2],z \in [0,1]\}$ denotes the product sutured manifold obtained by ``capping off'' one of the sutures of $A \times I$ as described in \cite[Defn. 2.8, Prop. 6.2]{GT08071432},
  \item $i: A \times I \longrightarrow D_A \times I$ denotes the inclusion map,
  \item $\hatbL := i(\bL) \subset D_A \times I$ denotes the image of $\bL$ under the inclusion, 
  \item $p(\bL) := \left(lk(\bL, B) \mod 2\right)$ denotes the ``mod 2 linking number'' of $\bL$ and $B$, which we will sometimes call the {\em parity} of $\bL$,
  \item $\Theta$ denotes the bigraded $\Ztwo$ vector space formally generated by $\theta^+,\theta^-$: \[\Theta := Span_{\Ztwo}(\theta^+,\theta^-),\] with bigrading given by: \begin{eqnarray*}
  deg(\theta^+) &=& (0,0),\\
  deg(\theta^-) &=& (-1,-1), \mbox{ and}
\end{eqnarray*}
  \item $\Thetahat$ denotes the singly-graded $\Ztwo$ vector space obtained by forgetting the first of the two gradings on $\Theta$.
%\begin{eqnarray*}({\bf A}(\theta^+),{\bf M}(\theta^+))&=& (0,0),\\
%                 ({\bf A}(\theta^-),{\bf M}(\theta^-))&=& (-1,-1),
%\end{eqnarray*}
\end{itemize} 

We will also periodically use the following abbreviated notation:
\begin{eqnarray*}
Y_{\bL} &:=& \boldSigma(A \times I, \bL),\,\, \mbox{and}\\
Y_{\hatbL} &:=& \boldSigma(D_A\times I, \hatbL).
\end{eqnarray*}

%In what follows, let 
%\begin{itemize}
%  \item $\bL \subset A \times I$ be a link, 
%  \item $\Theta = \langle \theta_+, \theta_-\rangle$ the $\Ztwo$ vector space defined as in the statement of Proposition \ref{prop:SFHFiltCpx}.
%\end{itemize} 

%Whenever $X$ is a bigraded vector space, we will denote by $X\{a,m\}$ the bigraded vector space whose bigrading has been shifted by $(a,m)$.  In particular, if ${\bf x} \in X$ is a bi-homogeneous element and ${\bf x}_{\{a,m\}} \in X\{a,m\}$ the corresponding element in the shifted vector space, then \[\left({\bf A}({\bf x}) + a\,\,,\,\,{\bf M}({\bf x}) + m\right)  = ({\bf A}({\bf x}_{\{a,m\}}),{\bf M}({\bf x}_{\{a,m\}})).\]

\begin{proposition} \label{prop:SFHFiltCpx}
$CFH(\boldSigma(D_A \times I, \widehat\bL))$ can be given the structure of a filtered chain complex whose associated graded complex has homology 

\[SFH(\boldSigma(A \times I, \bL)) \cong \left\{\begin{array}{ll}
                               \widehat{HFK}(\boldSigma(S^3, \bL), \widetilde{B}) \otimes \Theta & \mbox{ if $p(\bL) = 1$,}\\
			       \widehat{HFK}(\boldSigma(S^3, \bL), \widetilde{B}) & \mbox{ if $p(\bL)=0$,}\end{array}\right.\]

Moreover, the spectral sequence associated to this filtered complex converges to

 \[SFH(\boldSigma(D_A\times I,\hatbL)) \cong \left\{\begin{array}{ll}
	    \widehat{HF}(\boldSigma(S^3,\bL)) \otimes \Thetahat & \mbox{ if $p(\bL) = 1$, and}\\
            \widehat{HF}(\boldSigma(S^3,\bL))\left\{\frac{1}{2}\right\} \otimes \Thetahat & \mbox{ if $p(\bL) = 0$.}\end{array}\right.\]
\end{proposition}

\begin{remark}\label{remark:MandAgradings}
The bigrading referenced in the above proposition is the $({\bf A},{\bf M})$ bigrading associated to the null-homologous link $\widetilde{B} \subset \boldSigma(S^3, \bL)$, where ${\bf A}$ is the filtration (Alexander) grading on knot Floer homology defined in \cite{MR2065507}, \cite{GT0306378} and elaborated in \cite{GT0512286}, \cite{GT0607691} (see also \cite{GT064360}, \cite{GT07100359}), and ${\bf M}$ is the homological (Maslov) grading on Heegaard Floer homology defined in \cite{MR2222356}.

Note that an absolute (${\bf A} \in \Qrat$)--grading is only well-defined once a (homology class of) Seifert surface has been chosen for $\widetilde{B}$, and an absolute (${\bf M} \in \Qrat$)--grading is only well-defined in torsion Spin$^c$ structures of $\boldSigma(S^3, \bL)$.  (The proof of) Proposition \ref{prop:SFHFiltCpx} asserts the existence of a canonical identification between filtered chain complexes, preserving relative bigradings.\footnote{With respect to a fixed Seifert surface, relative ${\bf A}$--gradings lie in $\Zint$, while in a fixed Spin$^c$ structure, relative ${\bf M}$--gradings lie in some quotient of $\Zint$.  See \cite[Sec. 2]{MR2065507}, \cite[Sec.4]{MR2113019}.}  Therefore, $SFH(\boldSigma(A \times I, \bL))$ inherits an $({\bf A},{\bf M})$--bigrading (resp., $SFH(\boldSigma(D_A \times I,\hatbL))$ inherits an ${\bf M}$--grading), subject to these parameters.

%Note that when $\boldSigma(S^3, \bL)$ has positive first Betti number, the ${\bf A}$ grading (as an element of $\Zint$) depends on the choice of a Seifert surface for $\widetilde{B}$.  A canonical such choice is given by 

%(*********Maybe it's best not to dwell too much on this absolute bigrading.  Just note that the isomorphism alluded to in the proposition is canonical, and when absolute $({\bf A},{\bf M})$ gradings are defined, it respects these gradings.) When $\widehat{HF}(\boldSigma(S^3,\bL))$ is supported in torsion Spin$^c$ structures, as  An immediate corollary of Proposition \ref{prop:SFHFiltCpx} is that $SFH(\boldSigma(A\times I, \bL))$ comes equipped with a bigrading, $({\bf A},{\bf M}) \in \Qrat \times \Qrat$,  It is this bigrading which is referenced in the statement of the proposition.  
\end{remark}

\begin{proof}[Proof of Proposition \ref{prop:SFHFiltCpx}]
Note that $\boldSigma(A \times I,\bL)$ is the sutured manifold $(Y,\Gamma)$ where \[Y = \boldSigma(S^3,\bL) - N(\widetilde{B}) \mbox{ and } s(\Gamma) = \{\mu_1,\mu_1', \mu_2,\mu_2'\},\] where the $\mu_i$ represent meridians of $\widetilde{B}$ and the $\mu_i'$ represent oppositely-oriented meridians of $\widetilde{B}$.  To see this, observe that $A \times I$ is the standard sutured knot complement of $B \subset S^3$, hence its preimage under $\pi$ is the complement of $\widetilde{B} \subset \boldSigma(S^3,\bL)$.  It follows that each meridional (resp., oppositely-oriented meridional) suture lifts to two meridional (resp., oppositely-oriented meridional) sutures in $\boldSigma(A \times I, \bL)$.%  Furthermore, $\widetilde{B}$ is a $2$--component (resp., $1$--component) link if $p(\bL)=0$ (resp., $p(\bL)=1$).

Furthermore:
\begin{itemize}
\item when $p(\bL) = 0$, $\widetilde{B} = \widetilde{B}_1 \amalg \widetilde{B}_2$ is a $2$--component link.  $(Y,\Gamma)$ is then the standard sutured link complement (see \cite[Example 2.4]{MR2253454}), with a pair, $(\mu_i, \mu_i')$, of sutures on each boundary component, $-\partial(N(\widetilde{B}_i))$, and
\item when $p(\bL)=1$, $\widetilde{B}$ is a $1$--component link.  Hence, $(Y,\Gamma)$ is a sutured link complement with both pairs, $(\mu_i, \mu_i')$, of sutures on the single component of $\partial Y = -\partial(N(\widetilde{B}))$.  Note that this is a {\em non-standard} sutured link complement, since there is an extra pair of sutures on the single component of $-\partial(N(\widetilde{B}))$.
\end{itemize}

Similarly, we see that $\boldSigma(D_A \times I, \hatbL)$ is obtained from $\boldSigma(S^3,\bL)$ by deleting two $B^3$'s in $N(\widetilde{B})$.  Put differently, one obtains $\boldSigma(D_A \times I, \hatbL)$ from $\boldSigma(A \times I, \bL) = \boldSigma(S^3, \bL) - N(\widetilde{B})$ by filling in two $B^3$'s in $N(\widetilde{B})$.

Now, suppose that $(\Sigma, \boldalpha,\boldbeta)$ is a sutured Heegaard diagram for $\boldSigma(A \times I,\bL)$.  Then the topological observations above imply that one obtains a $4$--pointed Heegaard diagram \[(\Sigma',\boldalpha',\boldbeta',\mathbb{O}=(O_1,O_2),\mathbb{X} = (X_1,X_2))\] for $\boldSigma(S^3,\bL)$ compatible with the link $\widetilde{B}$ (in the sense of \cite{MR2065507}, \cite{GT0512286}, \cite{GT0607691}), by taking

\begin{itemize}
  \item \[\Sigma' := \Sigma \cup_{i=1,2} D_{O_i} \cup_{i=1,2} D_{X_i},\] where $D_{X_i}$ (resp., $D_{O_i}$) are disks containing a distinguished basepoint $X_i$ (resp., $O_i$) filling in the boundary components of $\Sigma$ corresponding to the meridians (resp., oppositely-oriented meridians) of $\widetilde{B}$,
  \item $\boldalpha'$ (resp., $\boldbeta'$), are the images of $\boldalpha$ (resp., $\boldbeta$) under the inclusion map $\Sigma \rightarrow \Sigma'$.
\end{itemize}

Similarly, by gluing in only the disks $D_{X_1}, D_{X_2}$ to $(\Sigma, \boldalpha,\boldbeta)$, one obtains a sutured Heegaard diagram, $(\hatSigma,\widehat{\boldalpha},\widehat{\boldbeta})$, for $\boldSigma(D_A \times I, \hatbL)$.

Let $\mathcal{C}'$ represent the Heegaard-Floer chain complex associated to $(\Sigma',\boldalpha',\boldbeta',\OO,\XX)$ with  

\begin{itemize}
  \item Generators: $\{{\bf x} \in \Torus_{\alpha'} \cap \Torus_{\beta'}\}$\\
  \item Differential: \[{\partial}'({\bf x}) = \sum_{{\bf y} \in \Torus_{\alpha'} \cap \Torus_{\beta'}}\sum_{\{\phi \in \pi_2({\bf x},{\bf y})\,|\,\mu(\phi) =1, n_{\OO}({\phi}) = 0\}} \widehat{M}(\phi)\cdot {\bf y},\]
\end{itemize}
where $n_{\OO}(\phi) := \sum_{i=1}^2 n_{O_i}(\phi)$.  Then the sutured Floer chain complex associated to $(\hatSigma,\widehat{\boldalpha}, \widehat{\boldbeta})$ is canonically isomorphic to $\mathcal{C}'$,\footnote{after fixing a generic choice of almost complex structure on $\hatSigma$, which extends uniquely to a generic choice of almost complex structure on $\Sigma'$} since the condition that $n_{\OO}(\phi) = 0$ for $\phi \subset \Sigma'$ is equivalent to the condition that $\phi \subset \hatSigma$ is never adjacent to $\partial\hatSigma$.  Hence, \[SFH(\boldSigma(D_A\times I,\bL)) \cong \left\{\begin{array}{ll}
	    \widehat{HF}(\boldSigma(S^3,\bL)) \otimes \Thetahat & \mbox{ if $p(\bL) = 1$, and}\\
            \widehat{HF}(\boldSigma(S^3,\bL))\left\{\frac{1}{2},\frac{1}{2}\right\} \otimes \Thetahat & \mbox{ if $p(\bL) = 0$,}\end{array}\right.\] as desired.  The extra $\left\{\frac{1}{2},\frac{1}{2}\right\}$ grading shift when $p(\bL) = 0$ arises due to the convention that $\widehat{HFK}(Y,L)$, for $L$ an $\ell$--component link, is defined to be the knot Floer homology of the ``knotification'' of $L$ in $Y \#^{\ell - 1} (S^1 \times S^2)$ (see \cite[Sec. 2.1]{MR2065507}).

Furthermore, since intersection positivity guarantees that $n_{\XX}(\phi) \geq 0$ for all holomorphic $\phi \in \pi_2({\bf x},{\bf y})$, we can endow $\mathcal{C}'$ with the structure of a $\Zint$--filtered complex by %on \[\Torus_{\widehat{\alpha}} \cap \Torus_{\widehat{\beta}} = \Torus_{\alpha'} \cap \Torus_{\alpha'} = \Torus_\alpha \cap \Torus_\beta\] 
%described in Remark \ref{remark:MandAgradings} endows $\mathcal{C}'$ with a $\Zint$--filtration.  More precisely, we can 
decomposing the differential as \[\partial' = \sum_{i=0}^\infty \partial'_i,\] where \[\partial'_i({\bf x}) := \sum_{{\bf y} \in \Torus_{\alpha'} \cap \Torus_{\beta'}}\sum_{\{\phi \in \pi_2({\bf x},{\bf y})\,|\,\mu(\phi) =1, n_{\OO}({\phi}) = 0,n_{\XX}({\phi}) = i\}} \widehat{M}(\phi)\cdot {\bf y}.\]

%$\phi \in \pi_2({\bf x},{\bf y})$ represent a  obtained by counting the intersection of any $\phi \in \pi_2({\bf x},{\bf y})$ complex, in the language of \cite[Sec. 2.4]{MR2065507}.  That the differential is a $\Qrat$--filtered morphism follows from intersection positivity and the fact that ${\bf A}({\bf x}) - {\bf A}({\bf y}) = n_\XX(\phi)-n_\OO(\phi)$ for $(\phi \subset \Sigma') \in \pi_2({\bf x},{\bf y})$.

The associated graded complex, $(\mathcal{C}',\partial'_0),$ is therefore chain isomorphic to the sutured Floer complex obtained from $(\Sigma, \boldalpha,\boldbeta)$, since counting $\phi \subset \Sigma'$ with $n_{\XX}(\phi) = n_{\OO}(\phi) = 0$ is equivalent to counting $\phi \subset \Sigma$ disjoint from $\partial\Sigma$.   Hence, $SFH(\boldSigma(A \times I,\bL))$ is the homology of the associated graded complex of the filtered complex $CFH(\boldSigma(D_A \times I, \hatbL))$.

This immediately implies (see \cite[Examples 2.3, 2.4]{MR2253454}, \cite[Prop. 9.2]{MR2253454}, and \cite[Thm. 1.1]{GT0512286}, combined with \cite[Lemma 2.12]{GT0610559}):
 
\[SFH(\boldSigma(A \times I, \bL)) \cong \left\{\begin{array}{ll}
                               \widehat{HFK}(\boldSigma(S^3, \bL), \widetilde{B}) \otimes \Theta & \mbox{ if $p = 1$,}\\
			       \widehat{HFK}(\boldSigma(S^3, \bL)), \widetilde{B}) & \mbox{ if $p=0$,}\end{array}\right.\]
as desired.

%Let $\mathcal{C}$ (resp., $\mathcal{C}'$, $\widehat{\mathcal{C}}$) represent the Heegaard-Floer chain complex associated to $(\Sigma,\boldalpha,\boldbeta)$ (resp., $(\Sigma',\boldalpha',\boldbeta',\OO,\XX)$, $(\hatSigma,\widehat{\boldalpha},\widehat{\boldbeta})$).   

%\begin{enumerate}
%  \item k
%\end{enumerate}
%is therefore canonically chain isomorphic to a chain complex, $\widetilde{C}$, associated to the $4$--pointed Heegaard diagram $(\Sigma',\boldalpha',\boldbeta')$ described in the proof of Proposition \ref{prop:SFHFiltCpx}, with

\end{proof}

%The above proposition tells us that $SFH(\boldSigma(A \times I, \bL))$ can be identified with (a variant of) the knot Floer homology of $\widetilde{B}$.  Since knot Floer homology is, in general, the associated graded complex of some filtered complex, it is natural to ask: Of what filtered complex is $SFH(\boldSigma(A \times I, \bL))$ the associated graded complex?  Proposition \ref{prop:SFHFiltCpx} provides the answer.

%In what follows, let \[\bL = \bL_{triv} \amalg \bL_{non} = (K_1 \amalg \ldots \amalg K_t) \amalg (K_{t+1} \amalg \ldots K_{t+n}) \subset A \times I\] be as in Section \ref{sec:KBSM} and $(\hatbL = i(\bL)) \subset D_A \times I$ the image of $\bL$ under the inclusion $i: A \times I \rightarrow D_A \times I$ as above.

We are now ready to define a sutured Floer-type ``functor''.

In what follows, let $\left(\bL = \bL_{triv} \amalg \bL_{non}\right) \subset A \times I$ be a resolved link with 
\begin{eqnarray*}
\bL_{triv} &=& K_1 \amalg \ldots \amalg K_t, \,\,\mbox{and}\\
\bL_{non} &=& K_{t+1} \amalg \ldots \amalg K_{t+n},
\end{eqnarray*}

and $\hatbL=i(\bL) \subset (D_A \times I)$ its image under the inclusion $i: A \times I \rightarrow D_A \times I$ as before.
\begin{definition} %For $\bL = K_1 \amalg \ldots \amalg K_{t+n} \subset A \times I$ a resolved link, %with 
%\begin{eqnarray*}
%\bL_{triv} = K_1 \amalg \ldots \amalg K_t, \mbox{ and}\\
%\bL_{non} = K_{t+1} \amalg \ldots \amalg K_{t+n}.
%\end{eqnarray*}
Let $\ZH(\bL)$ denote the $\Ztwo$--vector space formally generated by $[K_1],\ldots,[K_{t+n}]$: \[\ZH(\bL) := Span_{\Ztwo}([K_1],\ldots,[K_{t+n}]),\] endowed with the bigrading: 
\[deg([K_i]) :=\left\{\begin{array}{ll}
                        (0,-1) & \mbox{for } 1 \leq i \leq t, \mbox{and}\\
                        (-1,-1) & \mbox{else.}
                      \end{array}\right.\]

The space $\ZH(\bL)$ decomposes as $\ZH(\bL_{triv}) \oplus \ZH(\bL_{non})$, where
\begin{eqnarray*}
  \ZH(\bL_{triv}) &:=& Span_{\Ztwo}([K_1],\ldots,[K_{t}]),\\
  \ZH(\bL_{non}) &:=& Span_{\Ztwo}([K_{t+1}],\ldots,[K_{t+n}]).
\end{eqnarray*}
\end{definition}

\begin{definition} \label{defn:VHofL} Let $\VbarH(\bL)$ denote the exterior algebra \[\VbarH(\bL):=\Wedge^*\ZH(\bL),\] with its bigrading induced from $\ZH$: \[deg([K_{i_1}] \wedge \ldots \wedge [K_{i_k}]) := deg([K_{i_1}]) + \ldots + deg([K_{i_k}]),\] as in Definition \ref{defn:VofL}, and let $\VH(\bL)$ denote the following normalization of $\VbarH(\bL)$:\[\VH(\bL):=\VbarH(\bL)\left\{\frac{n-p(\bL)}{2},\frac{t+n-p(\bL)}{2}\right\}.\]
\end{definition}

%Denote by
%\begin{eqnarray*}
%Y_{\bL} &:=& \boldSigma(A \times I, \bL),\\
%Y_{\hatbL} &:=& \boldSigma(D_A \times I,\hatbL),\\
%\Wedge &:=& \Wedge^*(H^2(\#^{t+n} S^1\times S^2;\Ztwo)), \mbox{and}\\
%\Wedge_{triv} &:=& \Wedge^*(H^2(\#^t S^1 \times S^2;\Ztwo)).\\
%\end{eqnarray*}

Recall that $p(\bL):= lk(\bL,B) \,\mbox{mod}\, 2$.  In particular, when $\bL$ is a resolved link, \[p(\bL) = n \,\mbox{ mod }\, 2.\]  %$SFH(Y_{\bL};a)$ refers to the part of $SFH(Y_{\bL})$ in {\bf A} grading $a$.

\begin{proposition} \label{prop:SFHFunctor} Let $\bL \subset A \times I$ be a resolved link, as above.  Then \[SFH(Y_{\bL}) \cong \VH(\bL)\] as an $({\bf A},{\bf M})$--bigraded $\Ztwo$ vector space.

%Furthermore, 
%, is a free, rank $1$ module over $\Wedge$.  

%Furthermore,  $SFH(Y_{\bL})$, the homology of the associated graded complex of the $\Qrat$--filtered chain complex $CFH(Y_{\hatbL})$, decomposes as:

%\[SFH(Y_{\bL}) \cong \bigoplus_{a \in \Qrat} SFH\left(Y_{\bL};a\right),\] where 
%\[SFH(Y_{\bL};a) \cong \left\{\begin{array}{ll} \mbox{A free $\Lambda_{triv}$-module of rank $\binom{n}{k}$} & \mbox{ if $a \in \left\{\frac{n- p(\bL)}{2} - k\,\,|\,\, k \in \{0,\ldots, n\}\right\}$}\\
%0 & \mbox{ otherwise},\end{array}\right.\] where, for $a = \frac{n-p(\bL)}{2} - k$, all $\binom{n}{k}$ of the $\Wedge_{triv}$--module generators of $SFH(Y;a)$, are supported in {\bf M} grading \[\frac{(t+n) - p(\bL)}{2} - k.\]
\end{proposition}

\begin{proof}
%Let $(\hatbL = i(L)) \subset D_A \times I$ denote the image of $\bL$ under the inclusion $i:A \times I \rightarrow D_A \times I$, as before.  
Recall (see \cite[Sec. 4.2.5]{MR2113019}, \cite[Sec. 3.4]{GT08071432}) that if $(Y,\Gamma)$ is a balanced sutured manifold, then $SFH(Y,\Gamma)$ admits an action of $\Wedge^*(H_1^\circ(Y,\partial Y))$, where here (and in what follows), we use $H_*^\circ(-)$ to denote $H_*(-;\Zint)/\mbox{Tors}$. 

Since Proposition \ref{prop:SFHFiltCpx} tells us that $SFH(Y_{\bL})$ is the homology of the associated graded complex of $CFH(Y_{\hatbL})$, it will be convenient for us to begin by examining $SFH(Y_{\hatbL})$.

In particular, \cite[Lem. 5.16]{GT08071432} tells us that $SFH(Y_{\hatbL})$ is a free, rank $1$ module over 
\[\Wedge^*(H_1^\circ(Y_{\hatbL},\partial Y_{\hatbL}) \otimes_{\Zint} \Ztwo).\]

Furthermore, in the proof of \cite[Lem. 5.16]{GT08071432} we see that
\begin{eqnarray*}
 H_1^\circ(Y_{\hatbL},\partial Y_{\hatbL}) \otimes_{\Zint} \Ztwo &\cong& H^2(\#^{t+n} S^1 \times S^2;\Ztwo)\\
  &\cong& \ZH(\bL),
\end{eqnarray*} 

as ungraded vector spaces and hence \[\Wedge^*(H_1^\circ(Y_{\hatbL},\partial Y_{\hatbL}) \otimes_{\Zint} \Ztwo) \cong \VbarH(\bL)\]
%\begin{eqnarray*}
% \Wedge^*(H_1^\circ(Y_{\hatbL},\partial Y_{\hatbL}) \otimes_{\Z} \Ztwo) &\cong& \Wedge^*(H^2(\#^{t+n} S^1 \times S^2;\Ztwo))\\
%  &\cong& \VbarH(\bL)
%\end{eqnarray*} 
as ungraded algebras.  To see that they agree as bigraded algebras--and, hence, that $SFH(Y_\bL) \cong V_H(\bL)$ up to an overall bigrading shift--we need only show that under the identification: \[\phi:\ZH(\bL) \rightarrow H_1^\circ(Y_{\hatbL},\partial Y_{\hatbL}) \otimes_{\Zint} \Ztwo,\] elements of \[\ZH(\bL_{triv}) \leq \ZH(\bL) \,\,\,\mbox{(resp., of $\ZH(\bL_{non})\leq \ZH(\bL)$)}\] are sent to algebra generators of $({\bf A},{\bf M})$--degree $(0,-1)$, (resp., (-1,-1)).  %\[({\bf A},{\bf M}) = (0,-1) \,\,\,\mbox{(resp., an $({\bf A},{\bf M}) = (-1,-1)$)}\] bigrading shift on $SFH(Y_{\bL})$ 
%under the $H_1^\circ$--action.  

Since it is immediate, by definition of the $H_1^\circ$ action (see \cite[Sec. 4.2.5]{MR2113019}, \cite[Sec. 3.4]{GT08071432}), that the image of each element of $\ZH(\bL)$ has ${\bf M}$--degree $-1$, as desired, we are left to verify that the ${\bf A}$--degrees of elements of $\ZH(\bL)$ are as stated.

To see this, we will study the $H_1^\circ$ action on a chain complex, $CFH(Y_{\hatbL})$, associated to a particularly nice sutured Heegaard diagram for $Y_{\hatbL}$.  The generators of $CFH(Y_{\hatbL})$ will have well-defined absolute ${\bf A}$--gradings, and the homology of the associated graded complex will be $SFH(Y_\bL)$, as in Proposition \ref{prop:SFHFiltCpx}.

%$\clubsuit$  Note (7/8/2009):  The $H_1^\circ$ action is as claimed because (for algebra generators) it is defined by choosing a cochain which assigns $1$ to an element $\phi \in \Omega(T_\alpha,T_\beta)$ representing the correct element of $H_1^\circ$.  For the standard Heegaard diagram associated to connect sums of $S^1 \times S^2$, this assigns $1$ to either domain connecting the two intersection points associated to a compatible pair of $\alpha, \beta$ curves.   $A$ gradings are well defined on $SFH(Y_\bL)$ since if ${\bf x},{\bf y}$ are a pair of generators and $\phi,\phi' \in \pi_2({\bf x},{\bf y})$, then $\phi' = \phi + \Psi$, where $\Psi$ is a periodic domain, representing an element of $H_2(Y)$.  Since $K$ is nullhomologous, $#(K \cap [\Psi]) = 0$, by Poincare duality.
Accordingly, we proceed as in \cite[Sec. 5.3]{GT08071432} by choosing an isotopy class representative of ${\bL}$ satisfying $K_i \pitchfork A_{\frac{1}{2}} \neq \emptyset$ for all $i \in \{1, \ldots t+n\}$, where $A_{\frac{1}{2}} := A \times \left\{\frac{1}{2}\right\}.$

%To see that $SFH(Y_{\hatbL})$ is a free, rank $1$ module over $\Wedge$, we note that, in this situation, $\hatbL$ is an admissible, balanced, resolved tangle (ABR) in $D_A \times I$ which does not backtrack, in the sense of \cite[Defns. 5.1-5.5]{GT08071432}.  Therefore, $SFH(Y_{\hatbL})$ is a free, rank $1$ module over \[H_1^\circ(Y_{\hatbL},\partial Y_{\hatbL}) \otimes_\Z \Ztwo \cong \Wedge,\] by \cite[Lem. 5.16, 5.18]{GT08071432}.

%To see that the homology, $SFH(Y_{\bL})$, of the associated graded complex of the filtered complex, $CFH(Y_{\hatbL})$, decomposes as claimed, we will need to understand the ({\bf A},{\bf M}) bigrading.  $SFH(Y_{\bL})$ will then inherit an action of the subalgebra of $\Wedge$ preserving the {\bf A} grading.\footnote{This subalgebra of $\Wedge$, in turn, can be shown to agree with the part of $\Wedge^*(H_1^\circ(Y_{\bL},\partial Y_{\bL})\otimes_\Zint \Ztwo)$ which acts non-trivially on $SFH(Y_{\bL})$, in the sense of \cite[Sec. 3.4]{GT08071432}.  The corresponding actions also agree.} 

Let $\pi: Y_{\bL} \rightarrow A \times I$ be the double-branched covering projection, and $\pi_{\frac{1}{2}}: \pi^{-1}(A_{\frac{1}{2}}) \rightarrow A_{\frac{1}{2}}$ its restriction to the $\frac{1}{2}$ level.  Then $(\Sigma, \boldalpha,\boldbeta)$ is a sutured Heegaard diagram for $Y_{\bL}$, where 

\begin{itemize}
  \item $\Sigma = \pi_{\frac{1}{2}}^{-1}(A_{\frac{1}{2}})$
  \item $\boldalpha = \pi_{\frac{1}{2}}^{-1}(\cP(\mbox{cups}))$
  \item $\boldbeta = \pi_{\frac{1}{2}}^{-1}(\cP(\mbox{caps}))$
\end{itemize}

See Figures \ref{fig:SuturedHD} and \ref{fig:SuturedHD2} along with \cite[Sec. 5.3]{GT08071432} for more details.

\begin{figure}
\begin{center}
\resizebox{4.5in}{!}{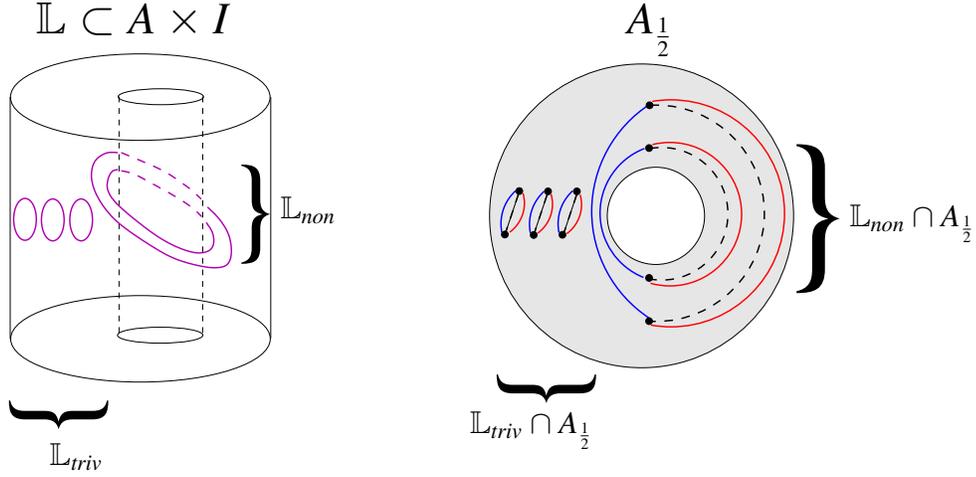}
\end{center}
\caption{The lefthand figure shows $(\bL = \bL_{triv} \amalg \,\,\bL_{non})$ standardly imbedded in $A \times I$.  The righthand figure shows the middle level, $A_{\frac{1}{2}}$.  The black dots represent the intersection points, $\bL \cap A_{\frac{1}{2}}$, the dotted lines are the branch cuts, and the red (resp., blue) arcs are the projections of the cups (resp., caps).}
\label{fig:SuturedHD}
\end{figure}

\begin{figure}
\begin{center}
\resizebox{4.5in}{!}{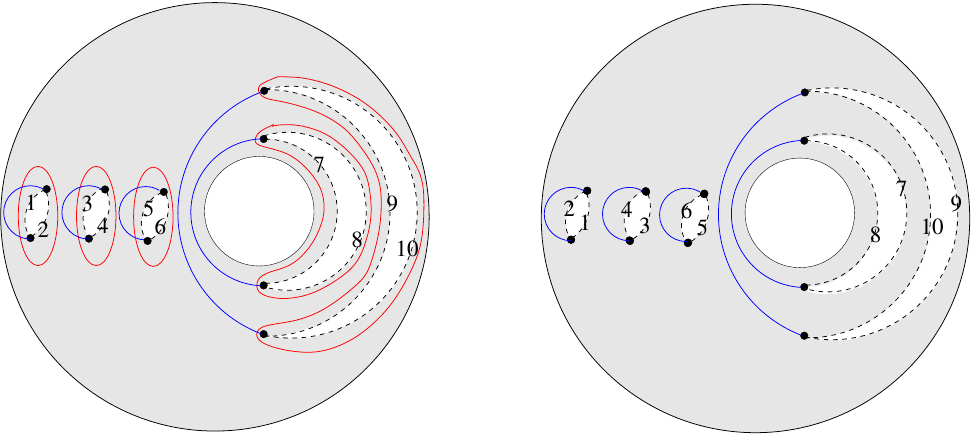}
\end{center}
\caption{A sutured Heegaard diagram, $(\Sigma, \boldalpha, \boldbeta)$, for $Y_{\bL}$, where $\bL$ is as in Figure \ref{fig:SuturedHD}.  The dotted lines are identified as labeled, and the colored circles are isotopic to the preimages of the colored arcs in Figure \ref{fig:SuturedHD}.}
\label{fig:SuturedHD2}
\end{figure}

Furthermore, as described in the proof of Proposition \ref{prop:SFHFiltCpx}, we obtain a sutured Heegaard diagram, $(\hatSigma,\widehat{\boldalpha},\widehat{\boldbeta})$, for $Y_{\hatbL} = \boldSigma(D_A \times I, \bL)$ by gluing in disks $D_{X_1}, D_{X_2}$ to two of the boundary components. Let 
\begin{eqnarray*}
  \widehat{\boldalpha} &=& \{\widehat{\alpha}_1, \ldots, \widehat{\alpha}_t, \widehat{\alpha}_{t+1}, \ldots \widehat{\alpha}_{t+n}\}, \,\,\, \mbox{and}\\
  \widehat{\boldbeta} &=& \{\widehat{\beta}_1, \ldots, \widehat{\beta}_t, \widehat{\beta}_{t+1}, \ldots, \widehat{\beta}_{t+n}\}
\end{eqnarray*}
be labeled as in Figure \ref{fig:Sigmahat}.  

\begin{figure}
\begin{center}
\resizebox{4.5in}{!}{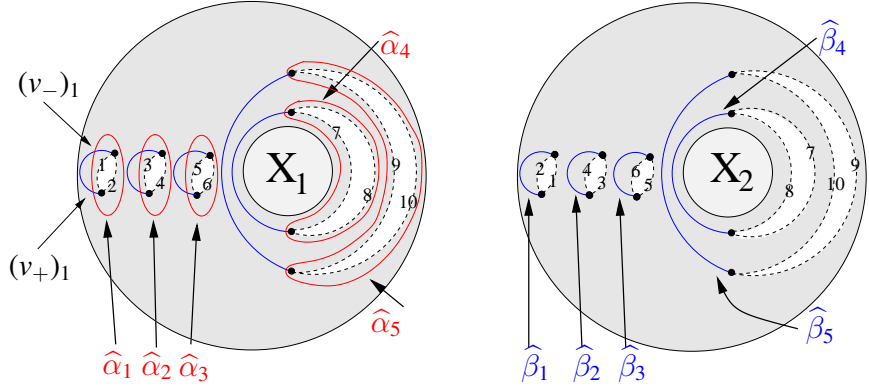}
\end{center}
\caption{A sutured Heegaard diagram for $Y_{\hatbL}$.  The $\widehat{\alpha}$, $\widehat{\beta}$ curves are numbered in increasing order from left to right (resp., from inside to outside) in the $\bL_{triv}$ (resp., the $\bL_{non}$) section of $\hatSigma$.}
\label{fig:Sigmahat}
\end{figure}

%\begin{figure}
%\begin{center}
%\input{Figures/Sigmahat.pstex_t}
%\end{center}
%\caption{A $4$--pointed Heegaard diagram for $(Y_{{\bL}_{non}},\widetilde{B})$, obtained by filling in the boundary components of $\Sigma$ in Figure \ref{fig:Sigmanon} with disks $D_{O_1},D_{O_2},D_{X_1},D_{X_2}$ (shaded lightly).  The darkly-shaded region is a domain, $\phi \in \pi_2({\bf x},{\bf y})$ for (for example) ${\bf x} = (v_+)_1 \otimes \ldots (v_+)_{n-1} \otimes (v_+)_n$ and ${\bf y} = (v_+)_1 \otimes \ldots (v_+)_{n-1} \otimes (v_-)_n.$}
%\label{fig:Sigmahat}
%\end{figure}

Since \[\alpha_i \cap \beta_j = \widehat{\alpha}_i \cap \widehat{\beta}_j = \left\{\begin{array}{ll}
                                  \{(v_+)_i, (v_-)_i\} & \mbox{if $i=j$}\\
				  \emptyset & \mbox{otherwise,}
				  \end{array}\right.\]

we obtain the following set of $\Ztwo$--module generators for both $CFH(Y_{\bL})$ and $CFH(Y_{\hatbL})$: \[\{(v_{\pm})_1 \otimes \cdots \otimes (v_{\pm})_{t+n}\}.\]

Furthermore, as in the proof of \cite[Lem. 5.18]{GT08071432}, $H_1^\circ(Y_{\hatbL},\partial Y_{\hatbL})$ has a basis whose elements, $\{[K_i]\}_{i=1}^{t+n}$, are in one-to-one correspondence with the components, $\{K_i\}_{i=1}^{t+n}$, of $\bL$,  %Now we use \cite[Lem. 5.18]{GT08071432}, which tells us that $\Wedge$ is formally generated as an exterior algebra by the components, $[K_1], \ldots, [K_{t+n}]$, of $\hatbL$.  
and the $H_1^\circ$-action on the filtered complex, $CFH(Y_{\hatbL})$ %obtained from $(\hatSigma,\widehat{\boldalpha},\widehat{\boldbeta})$ 
is extended linearly from the following action on generators:

Let $x_i \in \{(v_+)_i,(v_-)_i\}$ for all $i \in \{1, \ldots, t+n\}$, then 
\[[K_j] \cdot (x_1 \otimes \ldots \otimes x_{t+n}) := \left\{\begin{array}{ll}
                    x_1 \otimes \ldots \otimes (v_-)_j \otimes \ldots \otimes x_{t+n} & \mbox{if $x_j = (v_+)_j$, and}\\
		    0 & \mbox{if $x_j = (v_-)_j$.}\end{array}\right.\]

Recalling that if ${\bf x}, {\bf y} \in \Torus_{\widehat{\alpha}} \cap \Torus_{\widehat{\beta}}$ and $\phi \in \pi_2({\bf x},{\bf y})$, the relative ${\bf A}$ grading is given by: \[{\bf A}({\bf x}) - {\bf A}({\bf y}) = n_{\XX}(\phi)\]
%{\bf M}({\bf x}) - {\bf M}({\bf y}) &=& \mu(\phi),
%\end{eqnarray*}
where \[n_{\XX}(\phi):= \sum_{i=1}^2 n_{X_i}(\phi),\]% and $\mu(\phi)$ is the Maslov index of $\phi$, as usual.
there are two cases to consider.  %Since, by Equation \ref{eqn:RelAgrading}, multiplication by $[K_j]$ preserves the {\bf A} grading of a homogeneous element iff $j \in \{1, \ldots, t\}$, we see that the $\Wedge$ action descends to an action of $\Wedge_{triv}$ on the homology, $SFH(Y_{\bL})$, of the associated graded complex, as desired.

%To see that
%exterior algebra generators of $\Wedge^*(H_1^\circ(Y_{\hatbL},\partial Y_{\hatbL}))$
%To compute the absolute ${\bf A}$ and ${\bf M}$ gradings of these generators, it is most convenient to first compute relative ${\bf A}$ and ${\bf M}$ gradings, normalizing at the end.

%First recall that Proposition \ref{prop:SFHFiltCpx} gives us an expicit relationship between $SFH(Y_{{\bL}_{non}})$ and the knot Floer homology of $\widetilde{B}$.  

%In particular, let $(\Sigma, \boldalpha, \boldbeta)$ be a sutured Heegaard diagram for $Y_{{\bL}_{non}}$ as in Figure \ref{fig:Sigmanon} and let $(\Sigma',\boldalpha',\boldbeta',\OO,\XX)$ be the $4$--pointed Heegaard diagram obtained from $(\Sigma, \boldalpha,\boldbeta)$ by attaching disks $D_{O_1}, D_{O_2}, D_{X_1}, D_{X_2}$ to $\Sigma$ as specified in the proof of Proposition \ref{prop:SFHFiltCpx}.  See Figure \ref{fig:Sigmanonfilled}.  

%We consider two fundamental cases, which we use to compute relative {\bf A} and {\bf M} gradings between any pair of generators.

%\begin{remark}  The subalgebra of $\Wedge$ preserving the {\bf A} grading is strictly larger than $\Wedge_{triv}$ when $n\geq 2$.  However, it is only the subalgebra generated by $[K_i]$, $i \in \{1, \ldots, t\}$, described above, that acts nontrivially on 
%\end{remark}
%This is enough to conclude that $SFH(Y_{\bL};a)$, is as stated, up to an overall $({\bf A},{\bf M})$ bigrading shift.

{\flushleft {\bf Case 1: $K_j \in \bL_{triv}$}}

Let $x_i \in \{(v_+)_i,(v_-)_i\}$ for all $i = 1, \ldots, t+n$.  $[K_j]$ acts non-trivially on a generator, ${\bf x}$, iff it is of the form: \[{\bf x} = x_1 \otimes \ldots \otimes (v_+)_j \otimes \ldots \otimes x_{t+n},\] %and ${\bf y} = y_1 \otimes \ldots \otimes y_{t+n}$ are two generators in $\Torus_{\widehat{\alpha}} \cap \Torus_{\widehat{\beta}}$ 
in which case: \[[K_j]\cdot({\bf x}) = x_1 \otimes \ldots \otimes (v_-)_j \otimes \ldots \otimes x_{t+n}.\] 
%where, $x_i \in \{(v_+)_i,(v_-)_i\}$ for all $i\neq j$.
% satisfying the property that $\exists j \in \{1, \ldots, t\}$ such that \[\begin{array}{ll}
%        x_i = y_i & \mbox{ when $i \in \{1, \ldots, t+n\}, i\neq j$, and}\\
%	x_i = (v_+)_i \mbox{ and } y_i = (v_-)_i& \mbox{ when $i=j$.}
%\end{array}\]

%We claim that it follows that 
%\begin{eqnarray*}
%{\bf A}({\bf x}) - {\bf A}({\bf y}) &=& 0.\\
%{\bf M}({\bf x}) - {\bf M}({\bf y}) &=& 1.
%\end{eqnarray*}

Since ${\bf x}$ and $[K_j]\cdot({\bf x})$ are connected by a domain, $\phi$, as in Figure \ref{fig:Sigmahatgenstriv}, and $\phi$ satisfies $n_{\XX}(\phi) = 0$, we conclude that \[{\bf A}({\bf x})-{\bf A}([K_j]\cdot{\bf x}) = 0,\]  %Furthermore, $\phi$ is a bigon.  Using Lipshitz's formula for the Maslov index of a domain \cite[Cor. 4.3]{MR2240908}(see also \cite[Prop. 2.18]{GT08071432}) therefore tells us:
%\begin{eqnarray*}
%  {\bf M}({\bf x}) - {\bf M}({\bf y}) &=& \mu(\phi)\\
%                                      &=& [e(\phi)] + [n_{\bf x}(\phi) + n_{\bf y}(\phi)]\\
%                                      &=& \frac{1}{2} + \frac{1}{2}\\
%                                      &=& 1,
%\end{eqnarray*}
as desired.

\begin{figure}
\begin{center}
\resizebox{4.5in}{!}{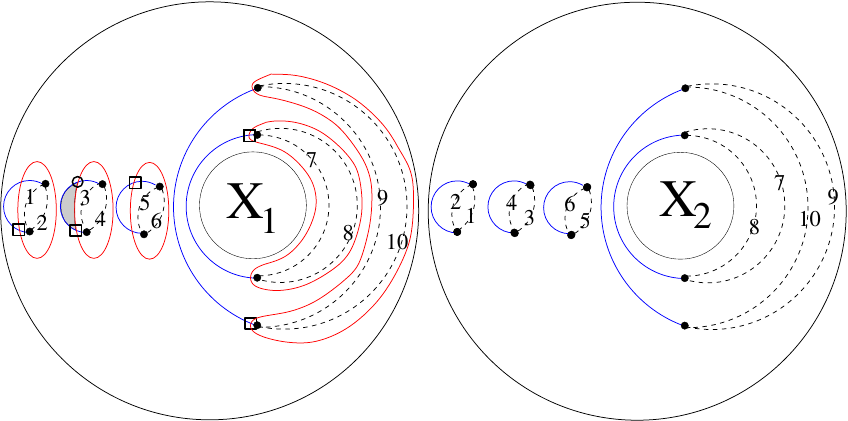}
\end{center}
\caption{The shaded region represents a domain, $\phi$, connecting two generators, {\bf x} and {\bf y}, in $(\hatSigma,\widehat{\boldalpha},\widehat{\boldbeta})$.  Here {\bf x} is represented by the square intersection point(s), while {\bf y} agrees with {\bf x} everywhere except at $\widehat{\alpha}_2 \cap \widehat{\beta}_2$, where {\bf y} is represented by the circular intersection point.}
\label{fig:Sigmahatgenstriv}
\end{figure}

{\flushleft {\bf Case 2: $K_j \in \bL_{non}$}}

As in Case 1, above, %Let $x_i \in \{(v_+)_i,(v_-)_i\}$ for all $i = 1, \ldots, t+n$.  
$[K_j]$ acts non-trivially on a generator, ${\bf x}$, iff it is of the form: \[{\bf x} = x_1 \otimes \ldots \otimes (v_+)_j \otimes \ldots \otimes x_{t+n},\] %and ${\bf y} = y_1 \otimes \ldots \otimes y_{t+n}$ are two generators in $\Torus_{\widehat{\alpha}} \cap \Torus_{\widehat{\beta}}$ 
in which case: \[[K_j]\cdot({\bf x}) = x_1 \otimes \ldots \otimes (v_-)_j \otimes \ldots \otimes x_{t+n}.\] 
%where, $x_i \in \{(v_+)_i,(v_-)_i\}$ for all $i\neq j$.
% satisfying the property that $\exists j \in \{1, \ldots, t\}$ such that \[\begin{array}{ll}
%        x_i = y_i & \mbox{ when $i \in \{1, \ldots, t+n\}, i\neq j$, and}\\
%	x_i = (v_+)_i \mbox{ and } y_i = (v_-)_i& \mbox{ when $i=j$.}
%\end{array}\]

%We claim that it follows that 
%\begin{eqnarray*}
%{\bf A}({\bf x}) - {\bf A}({\bf y}) &=& 0.\\
%{\bf M}({\bf x}) - {\bf M}({\bf y}) &=& 1.
%\end{eqnarray*}

When $K_j \in \bL_{non}$, then ${\bf x}$ and $[K_j]\cdot({\bf x})$ are connected by a domain, $\phi$, as in Figure \ref{fig:Sigmahatgensnon}. Since $\phi$ satisfies $n_{\XX}(\phi) = 1$, \[{\bf A}({\bf x})-{\bf A}([K_j]\cdot{\bf x}) = 1,\]  %Furthermore, $\phi$ is a bigon.  Using Lipshitz's formula for the Maslov index of a domain \cite[Cor. 4.3]{MR2240908}(see also \cite[Prop. 2.18]{GT08071432}) therefore tells us:
%\begin{eqnarray*}
%  {\bf M}({\bf x}) - {\bf M}({\bf y}) &=& \mu(\phi)\\
%                                      &=& [e(\phi)] + [n_{\bf x}(\phi) + n_{\bf y}(\phi)]\\
%                                      &=& \frac{1}{2} + \frac{1}{2}\\
%                                      &=& 1,
%\end{eqnarray*}
as desired.

%Now suppose that ${\bf x} = x_1 \otimes \ldots \otimes x_{t+n}$ and ${\bf y} = y_1 \otimes \ldots \otimes y_{t+n}$ are two generators satisfying the property that $\exists j \in \{t+1, \ldots, t+n\}$ such that \[\begin{array}{ll}
%        x_i = y_i & \mbox{ when $i \in \{1, \ldots, t+n\}, i\neq j$, and}\\
%	x_i = (v_+)_i \mbox{ and } y_i = (v_-)_i& \mbox{ when $i=j$.}
%\end{array}\]

%Then we claim 
%\begin{eqnarray*}
%{\bf A}({\bf x}) - {\bf A}({\bf y}) &=& 1, \mbox{ and}\\
%{\bf M}({\bf x}) - {\bf M}({\bf y}) &=& 1.
%\end{eqnarray*}  

%Then ${\bf x}$ and ${\bf y}$ are connected by a domain, $\phi$, as in Figure \ref{fig:Sigmahatgensnon} %which is topologically a $(j-t)$--punctured sphere satisfying 
%with $n_{\XX}(\phi) = 1$.  Hence, \[{\bf A}({\bf x})-{\bf A}({\bf y}) = 1,\] as desired.  %Using Lipshitz's formula for the Maslov index of a domain \cite[Cor. 4.3]{MR2240908}(see also \cite[Prop. 2.18]{GT08071432}), allows us to conclude that:
%\begin{eqnarray*}
%  {\bf M}({\bf x}) - {\bf M}({\bf y}) &=& \mu(\phi)\\
%                                      &=& [e(\phi)] + [n_{\bf x}(\phi) + n_{\bf y}(\phi)]\\
%                                      &=& \left[\frac{1}{2} - (n-1)\right] + \left[\frac{1}{2} + (n-1)\right]\\
%                                      &=& 1,
%\end{eqnarray*}

%as desired.

\begin{figure}
\begin{center}
\resizebox{4.5in}{!}{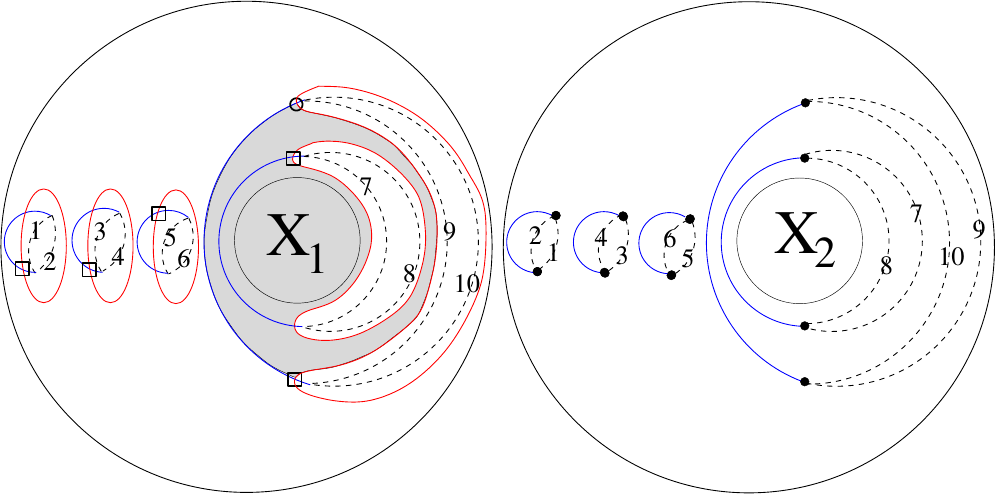}
\end{center}
\caption{The shaded region represents a domain, $\phi$, connecting two generators, {\bf x} and {\bf y}, in $(\hatSigma,\widehat{\boldalpha},\widehat{\boldbeta})$.  Here, {\bf x} is represented by the square intersection point(s), while {\bf y} agrees with {\bf x} everywhere except at $\widehat{\alpha}_5 \cap \widehat{\beta}_5$, where {\bf y} is represented by the circular intersection point.}
\label{fig:Sigmahatgensnon}
\end{figure}

%We can now compute the relative {\bf A} and {\bf M} gradings between any pair,\[{\bf x} = x_1 \otimes \ldots \otimes x_{t+n},{\bf y} = y_1 \otimes \ldots \otimes y_{t+n} \in \Torus_{\widehat{\alpha}} \cap \Torus_{\widehat{\beta}},\] of generators.  If $(x_+)$ (resp., $y_+$) denote the number of $x_i$ (resp., $y_i$) equal to $(v_i)_+$, then: 
%\begin{equation}\label{eqn:RelMgrading}
%{\bf M}({\bf x}) - {\bf M}({\bf y}) = x_+-y_+.
%\end{equation}

%Furthermore, if $(x_+)_{non}$ (resp., $(y_+)_{non}$) denotes the number of $x_i$ (resp., $y_i$) equal to $(v_i)_+$ for $i$ in the range $\{t+1, \ldots, t+n\}$, then: 
%\begin{equation}\label{eqn:RelAgrading}
%{\bf A}({\bf x}) - {\bf A}({\bf y}) = (x_+)_{non} - (y_+)_{non}.
%\end{equation}

We conclude that, with respect to this choice of sutured Heegaard diagram, \[CFH(Y_\bL) \cong \VH(\bL)\] as relatively $({\bf A},{\bf M})$ bigraded $\Ztwo$ modules, and, hence, that $SFH(Y_\bL) \cong \VH(\bL)$ as relatively $({\bf A},{\bf M})$ bigraded modules, since: \[2^{t+n} = rk(CFH(Y_\bL)) \geq rk(SFH(Y_\bL)) \geq rk(SFH(Y_{\hatbL})) = 2^{t+n}.\] Note that the second inequality above arises because $CFH(Y_{\bL})$ can be identified as the associated graded complex of the filtered complex $CFH(Y_{\hatbL})$, as in Proposition \ref{prop:SFHFiltCpx}.%, hence $SFH(Y_\bL) \cong CFH(Y_\bL) \cong \VbarH(\bL)$. 

To see that the absolute {\bf A} and {\bf M} gradings on $SFH(Y_\bL)$ are as stated, it suffices to show that the extremal generator, \[{\bf v}_+ := (v_+)_1 \otimes \ldots \otimes (v_+)_{t+n},\] has bigrading \[({\bf A},{\bf M}) = \left(\frac{n-p(\bL)}{2},\frac{t+n-p(\bL)}{2}\right).\]  

To see that ${\bf A}({\bf v}_+)$ is as claimed, use Proposition \ref{prop:SFHFiltCpx}, combined with the symmetry properties of knot Floer homology with extra basepoint pairs \cite[Sec. 3.5]{MR2065507}, \cite[Sec. 8]{GT0512286}, \cite[Lem. 2.12]{GT0610559}, which tell us that for a multiply-pointed Heegaard diagram associated to an $\ell$ component link $K = \amalg_{i=1}^\ell K_i$ for which there are $n_i$ $(\OO,\XX)$ basepoint pairs associated to the link component $K_i$, we have: \[\left(\widehat{HFK}(Y,K) \otimes \Theta^{\otimes s}\right)_{(a,m)} \cong \left(\widehat{HFK}(Y,K) \otimes \Theta^{\otimes s}\right)_{(-a-s,m-2a-s)}\]
where $s: = \sum_{i=1}^\ell (n_i -1)$, $\Theta$ is the bigraded $\Ztwo$ vector space generated by $\theta^+,\theta^-$ satisfying 
\begin{eqnarray*}({\bf A}(\theta^+),{\bf M}(\theta^+))&=& (0,0),\\
                 ({\bf A}(\theta^-),{\bf M}(\theta^-))&=& (-1,-1),
\end{eqnarray*} 
and $(\widehat{HFK}(Y,K) \otimes \Theta^{\otimes s})_{(a,m)}$ denotes the summand in $({\bf A},{\bf M})$ bigrading $(a,m)$.

Since the above symmetry in the case of interest ($Y = \boldSigma(S^3,\bL)$, $K=\widetilde{B}$) identifies ${\bf v}_+$ with ${\bf v}_- = (v_+)_1 \otimes \ldots (v_+)_t \otimes (v_-)_{t+1} \otimes \ldots \otimes (v_-)_{t+n}$, we obtain 
\begin{eqnarray*}
{\bf A}({\bf v}_+) - {\bf A}({\bf v}_-) &=& {\bf A}({\bf v}_+) - \left(-{\bf A}({\bf v}_+) - s\right)\\
                                      n &=& 2{\bf A}({\bf v}_+) + s\\
                          \frac{n-s}{2} &=& {\bf A}({\bf v}_+)\\
                          \frac{n-p(\bL)}{2} &=& {\bf A}({\bf v}_+),
\end{eqnarray*}

as desired.

To see that \[{\bf M}({\bf v}_+) = \frac{t+n-p(\bL)}{2},\]  we use the fact that \[\boldSigma(S^3,\bL) = \#^{(t+n)-1}(S^1 \times S^2),\] combined with the fact that, by convention, the knot Floer homology of a multi-component link in $Y$ comes from the ``knotification'' of the link \cite[Sec. 2.1]{MR2065507} in $Y \# (\ell -1) S^1 \times S^2$.  Using the Heegaard Floer K{\"u}nneth principle (\cite[Thm. 1.5]{MR2113020}) and the fact (cf. \cite[Subsec. 9.1]{MR2113019}) that $\widehat{HF}(S^1 \times S^2)$ is the free $\Ztwo$--module on generators $\{{\bf w}_+,{\bf w}_-\}$ with \[{\bf M}({\bf w}_\pm) = \pm \frac{1}{2}\] implies that \[{\bf M}({\bf v}_+) = \frac{t+n - p(\bL)}{2},\] as desired. $\qed$
\end{proof}

The following proposition establishes the equivalence of the Khovanov and sutured Floer ``functors'' on resolved links, $\bL \subset A \times I$.  In what follows, 
\begin{itemize}
  \item $V(\bL)$ (resp., $V_H(\bL)$) denotes the $(f,q)$--bigraded (resp., $({\bf A},{\bf M})$--bigraded) vector space of Definition \ref{defn:VofL} (resp., Definition \ref{defn:VHofL}).
  \item For resolved links $\bL',\bL'' \subset A \times I$ connected by a saddle cobordism $S$, $G\cF_S$ denotes the map \[SFH(Y_{\bL'}) \rightarrow SFH(Y_{\bL''})\] induced by performing surgery on the corresponding knot in $Y_{\bL'}$ as in \cite[Sec. 4]{GT08071432}.  %Note that if $\hatbL',\hatbL'' \subset D_A \times I$ are the images of $\bL',\bL''$ under the inclusion $i:A \times I \rightarrow D_A \times I$, then  this map is the (induced map on homology of the) ${\bf A}$--grading-preserving part of the map \[CFH(Y_{\hatbL'}) \rightarrow CFH(Y_{\hatbL''})\]  induced by the same surgery on the image of the same knot under the inclusion $Y_{\bL'} \rightarrow Y_{\hatbL'}$.  This follows from the fact that a sutured multi-diagram, $(\hatSigma,\widehat{\boldalpha},\widehat{\boldeta}',\widehat{\boldeta}'')$, for $Y_{\hatbL'}$ and $Y_{\hatbL''}$ can be obtained from a sutured multi-diagram for $Y_{\bL'}$ and $Y_{\bL''}$ by filling in two of the boundary components with disks labeled with basepoints $\XX = X_1, X_2$.  Domains $\phi \subset \Sigma$ representing holomorphic triangles are then in one-to-one correspondence with domains $\widehat{\phi} \subset \hatSigma$ representing ${\bf A}$--grading-preserving triangles, since ${\bf A}$--grading-preserving domains in $\hatSigma$ satisfy $n_{\XX}(\widehat{\phi}) = 0$ and hence correspond to domains in $\Sigma$ which miss the boundary.  For more details, see the proof of Proposition \ref{prop:SFHFiltCpx}.
  \item $G\cV_S$ denotes the $f$--grading-preserving part of the map induced by the saddle-cobordism $S$, as in Definitions \ref{defn:KhMergeMap} and \ref{defn:KhSplitMap}.
%To a saddle cobordism between two resolved link projections, $\cP(\bL'), \cP(\bL'') \subset A$ we associate the induced map \[SFH(Y_{\bL'}) \rightarrow SFH(Y_{\bL''})\] on sutured Floer homology obtained by counting triangles in the associated sutured multi-diagram (cf. \cite[Sec. 3 \& 4]{GT08071432}).  More precisely, if $\cP({\bL}')$ and $\cP({\bL}'')$ are related by a saddle cobordism, then $Y_{\bL''}$ can be obtained from $Y_{\bL'}$ by means of a single surgery on an imbedded knot in $Y_{\bL'}$.  After constructing a sutured Heegaard multi-diagram $(\Sigma,\boldalpha,\boldeta',\boldeta'')$ subordinate to this knot as in~\cite[Sec. 4]{GT08071432}, the map \[SFH(Y_{\alpha,\eta'}) \rightarrow SFH(Y_{\alpha,\eta''})\] is given as in \cite[Sec. 3.5]{GT08071432}.
%\item Note that

\end{itemize}
\begin{proposition} \label{prop:Equiv}For each resolved link \[\bL = (\bL_{triv} \amalg \bL_{non}) \subset A \times I\] there is a canonical isomorphism \[\Psi_{\bL}: V(\bL) \rightarrow SFH(Y_{\bL}).\]%and associated induced link \[\hatbL \subset D \times I,\] there is a canonical filtered chain homotopy equivalence: \[\psi_{\hatbL}: CV(\hatbL) \stackrel{\cong}{\longrightarrow} CFH(Y_{\hatbL})\]
%sending $(f,q)$ bigradings to \[({\bf A},{\bf M}) = \left(\frac{f - p(\bL)}{2},\frac{q - p(\bL)}{2}\right)\] bigradings which, in particular, induces an isomorphism on the homology of the associated graded complexes there is a canonical isomorphism: \[\psi_{\bL}: V(\bL) \stackrel{\cong}{\longrightarrow} SFH(Y_{\bL}).\]  
Furthermore, this isomorphism is natural in the sense that if $\bL', \bL'' \subset A \times I$ are resolved links and $S:{\bL'}\rightarrow {\bL''}$
is a saddle cobordism, the following
diagram commutes:
$$
\begin{CD}
V({\bL'})@>{G\cV_S}>>V({\bL''})\\
@V{\psi_{\bL'}}VV @VV{\psi_{\bL''}}V \\
SFH(Y_{\bL'})@>{G\cF_S}>> 
SFH(Y_{\bL''})
\end{CD}
$$
\end{proposition}

\begin{proof}
Let $\bL \subset A \times I$ be any resolved link.  Proposition \ref{prop:SFHFunctor} tells us that $SFH(Y_\bL) \cong V_H(\bL)$, so the canonical isomorphism, $\Psi_\bL: V(\bL) \rightarrow SFH(Y_\bL)$ for any resolved $\bL \subset A \times I$ is the one induced by the canonical isomorphism $Z(\bL) \rightarrow Z_H(\bL)$.

Now denote by $\hatbL$ the image of $\bL$ under the inclusion $i: A \times I \longrightarrow D_A \times I$.  Then $\hatbL$ is an admissible, $0$--balanced, resolved tangle in $D_A \times I$, in the sense of \cite[Defn. 5.1-5.5]{GT08071432}, so \cite[Prop. 5.17]{GT08071432} gives us an isomorphism \[\psi_{\hatbL}: V(\hatbL) \longrightarrow SFH(Y_{\hatbL})\]

%$\bL'$, $\bL''$ are resolved links in $A \times I$ connected by an saddle cobordism $S$.  Denote by $\hatbL'$ (resp., $\hatbL''$) the image of $\bL'$ (resp., $\bL''$) under the inclusion $i: A \times I \longrightarrow D_A \times I$.  

%$\hatbL', \hatbL''$ are ABR's in $D_A \times I$, in the sense of \cite[Defn. 5.1-5.5]{GT08071432}, so \cite[Prop. 5.17]{GT08071432} gives us isomorphisms 
%\begin{eqnarray*}
%\psi_{\hatbL'}: V(\hatbL') &\longrightarrow& SFH(Y_{\hatbL'})\\
%\psi_{\hatbL''}: V(\hatbL'') &\longrightarrow& SFH(Y_{\hatbL''})\\
%\end{eqnarray*}

such that if $\bL'$, $\bL''$ are resolved links in $A \times I$ connected by a saddle cobordism $S$, and $\hatbL'$, $\hatbL''$ their induced images in $D_A \times I$, then the following diagram commutes:

\[\begin{CD}
V({\hatbL'})@>{\cV_S}>>V({\hatbL''})\\
@V{\psi_{\hatbL'}}VV @VV{\psi_{\hatbL''}}V \\
SFH(Y_{\hatbL'})@>{\cF_S}>> 
SFH(Y_{\hatbL''})
\end{CD}\]

Furthermore, for any resolved $\bL \subset A \times I$, $V(\hatbL) \cong V(\bL)$ (resp., $SFH(Y_{\hatbL}) \cong SFH(Y_{\bL})$) as $\Ztwo$ vector spaces.  In fact, %up to an overall linear transformation on the grading, 
$V(\hatbL)$ is obtained from $V(\bL)$ (resp., $SFH(Y_{\hatbL})$ is obtained from $SFH(Y_{\bL})$) by forgetting the $f$--grading (resp., ${\bf A}$--grading), and the isomorphism $\psi_{\hatbL}$ agrees with $\Psi_{\bL}$ subject to this forgetful map.  Lastly, given any saddle cobordism $S$ connecting a resolved link $\bL' \subset A \times I$ to $\bL'' \subset A \times I$, $G\cV_S$ (resp., $G\cF_S$) is the $f$--grading-preserving (resp., ${\bf A}$--grading-preserving) part of the map $\cV_S$ (resp., $\cF_S$) corresponding to the saddle cobordism connecting $\widehat{\bL'}$ to $\widehat{\bL''}$.

%In fact, we can alternatively view this isomorphism of $\Wedge$-modules as a filtered chain map, by pulling back ({\bf A},{\bf M}) bigradings on the $\Ztwo$--generators of $V(\hatbL)$ via the isomorphism, $\Psi_{\hatbL}$.  More precisely, if ${\bf x} \in V(\hatbL)$, then we define \[({\bf A}({\bf x}), {\bf M}({\bf x})) := ({\bf A}(\psi_{\hatbL}({\bf x})), {\bf M}(\psi_{\hatbL}({\bf x}))).\]

%(*********Now what I want to say is that $GV(\bL)$ can be identified with the associated graded complex of this filtered chain complex, by the proper identification of the $(f,q)$ bigradings with the $({\bf A},{\bf M})$ bigradings.*******)  
% $\Psi_{\hatbL}$ (there, $\Psi_T$), which induces an isomorphism on homology.  If we can show that the gradings correspond as claimed, $\Psi_{\hatbL}$ will then be a filtered chain map.  By a direct examination of the Heegaard diagram for $Y_{\bL}$, it can be seen that the associated graded chain complexes, $V(\bL)$ and $CFH(Y_{\bL})$ are chain isomorphic, which then implies (cf. \cite[]{}) that $\Psi_{\hatbL}$ is, in fact, a filtered chain isomorphism, as claimed.  Furthermore, naturality of the maps $\Psi_{\bL}$ as described above will follow from naturality of the maps $\Psi_{\hatbL}$ (where we forget about the filtration), combined with the stated grading correspondence.

Thus, to show that
\[\begin{CD}
V({\bL'})@>{G\cV_S}>>V({\bL''})\\
@V{\psi_{\bL'}}VV @VV{\psi_{\bL''}}V \\
SFH(Y_{\bL'})@>{G\cF_S}>> 
SFH(Y_{\bL''})
\end{CD}\]

is a commutative diagram it suffices to show that for any two $f$--homogeneous generators ${\bf x} \in V(\bL')$ and ${\bf y} \in V(\bL'')$, \[f({\bf x}) = f({\bf y}) \,\,\Leftrightarrow\,\, {\bf A}(\Psi_{\bL'}({\bf x})) = {\bf A}(\Psi_{\bL''}({\bf y})).\]

But Definitions \ref{defn:VofL} and \ref{defn:VHofL} %and Proposition \ref{prop:SFHFunctor} 
imply that 
\begin{eqnarray*}
  {\bf A}(\Psi_{\bL'}({\bf x})) &=& \frac{f({\bf x}) - p(\bL')}{2}, \mbox{ and}\\
  {\bf A}(\Psi_{\bL''}({\bf y})) &=& \frac{f({\bf y}) - p(\bL'')}{2},
\end{eqnarray*}

so $f$ gradings agree iff the corresponding ${\bf A}$ gradings agree.  The result follows.%\[f({\bf x}) = f({\bf y}) \,\,\Leftrightarrow\,\, {\bf A}(\Psi_{\bL'}({\bf x})) = {\bf A}(\Psi_{\bL''}({\bf y})),\] as desired.

\end{proof}

We are now ready to prove Theorem \ref{thm:Reinterp}.

\begin{proof}[Proof of Theorem \ref{thm:Reinterp}] We proceed exactly as in the proof of \cite[Thm. 5.19]{GT08071432}.  Recall (see Definition \ref{defn:GrKhChainCpx} and Remark \ref{rmk:FrameDep})
%Remark \ref{rmk:fgrading})
that if $\bL \subset A \times I$ is a link and $\overline{\bL}$ is its mirror, then $Kh^*(\overline{\bL}):=V(\bL)$, where $V(\bL)$ is the homology of the cubical complex constructed by 
\begin{enumerate}
  \item numbering the crossings of a projection, $\cP(\bL)$, of $\bL$,
  \item forming the associated cube of resolutions whose vertices are the resolved links, $\cP_\cI(\bL)$, described in Section \ref{sec:KBSM}, arranged so that two vertices $\cI, \cI'$ are connected by an edge iff $\cP_\cI(\bL)$, $\cP_\cI'(\bL)$ are connected by a saddle cobordism, $S$, 
  \item assigning to each vertex of the cube the vector space $V(\cP_\cI(\bL))$ and to each edge of the cube the differential map $G\cV_S$.
\end{enumerate}

Similarly, using the link surgeries spectral sequence for sutured Floer homology (see \cite[Sec. 4]{GT08071432}) associated to the link in $Y_\bL$ formed by taking the preimage of simple arcs at each crossing of the projection $\cP(\bL)$ as in \cite[Fig. 7]{GT08071432}, we see that $SFH(Y_\bL)$ is the homology of the filtered cubical complex constructed by proceeding as in steps (1) and (2) above, but in step (3) assigning the vector space \[SFH(Y_{\cP_\cI(\bL)}) \cong V_H(\cP_\cI(\bL))\] to each vertex of the cube, assigning the map $G\cF_S$ to each edge of the cube, and assigning higher polygon maps as in \cite[Sec. 4]{GT08071432} to higher--dimensional faces of the cube.

The $E^2$ page of this link surgeries spectral sequence is the homology of the chain complex whose underlying vector space is \[\bigoplus_{\cI \in \{0,1\}^\ell} SFH\left(Y_{\cP_\cI(\bL)}\right)\] and whose differential consists solely of those maps assigned to the edges of the cube.  Proposition \ref{prop:Equiv} yields a canonical isomorphism between this chain complex and the chain complex whose homology is $V(\bL)\cong Kh^*(\overline{\bL})$.  The result follows.
%\begin{itemize}
%  \item underlying vector space is the direct sum of the vector spaces, $V_H(\cP_\cI(\bL))$, assigned to the vertices of the cube,
%  \item differential counts polygons in the sutured multi-diagram associated to the cube of resolutions.
%\end{itemize}

\end{proof}

\section{A relationship between annular links and balanced tangles}

We close by establishing the relationship between Roberts' spectral sequence, \cite[Prop. 1.1]{GT07060741}, for annular links (or, more precisely, its reinterpretation in Theorem \ref{thm:Reinterp}) and the spectral sequence for balanced tangles introduced in \cite[Prop. 5.20]{GT08071432}.  The following is a direct consequence of \cite[Thm. 4.5]{SurfDecomp}, which generalizes Juh{\'a}sz's surface decomposition theorem, \cite[Thm. 1.3]{MR2390347}, to filtered complexes arising in the context of link surgeries spectral sequences.  See \cite{SurfDecomp} for a more complete discussion of related naturality results. 
 
\begin{theorem} \label{thm:Cutting} Let $\bL \subset A \times I$ be an isotopy class representative of an annular link admitting a projection, $\cP(\bL)$, and let $\lambda \subset A$ be a properly imbedded oriented arc representing a nontrivial element of $H_1(A,\partial A)$ such that $\lambda$ intersects $\cP(\bL)$ transversely.   % $\lambda \times I$ intersects $\bL$ transversely.  
Let $T \subset D \times I$ be the balanced tangle in $D \times I$ obtained by decomposing along the surface, $\lambda \times I$, endowed with the product orientation.

Then the spectral sequence \[Kh^*(\overline{T}) \rightarrow SFH(\boldSigma(D \times I, T))\] is a direct summand of the spectral sequence \[Kh^*(\overline{\bL}) \rightarrow SFH(\boldSigma(A \times I, \bL)).\]

Furthermore, the direct summand is trivial if there exists some $\bL' \subset A \times I$ isotopic to $\bL$ satisfying \[|(\lambda \times I) \pitchfork \bL'| \lneq |(\lambda \times I) \pitchfork \bL|\]
\end{theorem}

\begin{proof}[Proof of Theorem \ref{thm:Cutting}] 
Let 
\begin{itemize}
  \item $\vec{p} :=  (\lambda \times I) \pitchfork \bL$ denote the collection of transverse intersection points of $\bL$ with $\lambda \times I$, 
  \item $S := \boldSigma(\lambda \times I, \vec{p})$ denote the preimage of $\lambda \times I$ in $\boldSigma(A \times I, \bL)$, and 
  \item $L$ denote the link in $\boldSigma(A \times I, \bL)$ obtained as the preimage of trivial arcs at each crossing of $\bL$ as in \cite[Fig. 7]{GT08071432}.
\end{itemize}

Then $S, L \subset \boldSigma(A \times I, \bL)$ satisfy the assumptions of \cite[Thm. 4.5]{SurfDecomp}.  Furthermore, the decomposition of $\boldSigma(A \times I, \bL)$ along $S$ is precisely the sutured manifold $\boldSigma(D \times I, T)$, and $L' \subset \boldSigma(D \times I, T)$, the image of $L$ in the decomposed manifold, is precisely the link in $\boldSigma(D \times I, T)$ obtained as the preimage of the trivial arcs associated to the crossings of $\cP(T) \subset D$. 

Now \cite[Thm. 4.5]{SurfDecomp} implies that the filtered complex, $X(L')$, associated to $L'$ is a direct summand of the filtered complex, $X(L)$, associated to $L$.  But $X(L')$ is precisely the filtered complex whose associated spectral sequence has $E^2$ term $Kh^*(\overline{T})$ and converges to $SFH(\boldSigma(D \times I, T)),$ while the spectral sequence associated to $X(L)$ has $E^2$ term $Kh^*(\overline{\bL})$ and converges to $SFH(\boldSigma(A \times I, \bL))$, as desired.

%As in Section \ref{sec:SFH}, one obtains a natural framed link, $L \subset \boldSigma({\bL})$, in the strongly-balanced sutured manifold, $\boldSigma(\bL)$, from the data of $\bL \subset A \times I$.  %Therefore, any full sutured multi-diagram compatible with $L$ induces a link surgeries spectral sequence with $E^2$ term $Kh(\bL)$ and $E^\infty$ term $SFH(\boldSigma(\bL))$.

%Therefore, if $L'$ denotes the induced image of $L$ in the sutured manifold obtained by decomposing along $S$, there exist $(0,1)$ sutured multi-diagrams for $L$ and $L'$ such that $(X^{(0,1)})'$ is a filtered subcomplex of $X^{(0,1)}$.

%To see that this implies the existence of a spectral sequence from $Kh(\bL)$ to $SFH(\boldSigma(\bL))$ containing a direct summand that is a spectral sequence from $Kh(T)$ to $SFH(\boldSigma(T))$, one need only note that:
%\begin{itemize}
%  \item The sutured manifold obtained by decomposing along $S$ is precisely $\boldSigma(T)$, where $T \subset D \times I$ is the balanced tangle obtained by cutting $\bL \subset A \times I$ along $\lambda \times I$.
%  \item The spectral sequence associated to the filtered complex, $X^{(0,1)}$ (resp., $(X^{(0,1)})'$) has $E^2$ term $Kh(\bL)$ (resp., $Kh(T)$) and $E^\infty$ term $SFH(\boldSigma(\bL))$ (resp., $SFH(\boldSigma(T))$).
%\end{itemize}

Furthermore, the direct summand, $X(L') \leq X(L)$, is trivial if the intersection of $\lambda \times I$ with $\bL$ is non-minimal, since we know that if $|\bL \pitchfork (\lambda \times I)|$ is non-minimal among isotopy class representatives of $\bL$, then all resolutions of $T$ backtrack, in the sense of \cite[Defn. 5.8]{GT08071432}, so the associated chain complex for $Kh(T)$ has no generators, hence the entire filtered complex, $X(L')$, and its associated spectral sequence are trivial.
\end{proof}

\bibliography{KBSM_Annulus}
\end{document}